\documentclass{article}
\usepackage{amsmath}
\usepackage{amsfonts}
\usepackage{amssymb}
\usepackage{amsthm}
\usepackage{array}

\usepackage{parskip} 

\usepackage{graphicx}
\usepackage[labelfont=bf]{caption}
\usepackage[top=1in, bottom=1in, left=1in, right=1in]{geometry}
\usepackage{hyperref}
\usepackage{enumerate}

\usepackage{listings} 
\usepackage{color} 
\usepackage[usenames,dvipsnames]{xcolor}
\usepackage{caption} 
\usepackage{framed} 

\definecolor{dkgreen}{rgb}{0,0.6,0}
\definecolor{lgreen}{rgb}{0.25,1,0}
\definecolor{purple}{rgb}{0.35,0.02,0.48}

\usepackage[numbered,framed]{matlab-prettifier}
\lstdefinestyle{mat}{
  frame=single,
  language=Matlab,
  showstringspaces=false,
}
\usepackage[utf8]{inputenc}
\usepackage[T1]{fontenc}
\usepackage{lmodern}

\newcommand{\reals}{\mathbb{R}}

\newtheoremstyle{mystuff}{}{}{\itshape}{}{\bfseries}{:}{.5em}{}
\theoremstyle{mystuff}

\newtheorem*{definition*}{Definition}

\newtheorem*{theorem*}{Theorem}

\newtheorem*{lemma*}{Lemma}
\newtheorem*{proposition*}{Proposition}

\newtheoremstyle{myexample}{}{}{}{}{\bfseries}{:}{.5em}{}
\theoremstyle{myexample}
\newtheorem*{example*}{Example}

  

\newtheoremstyle{named}{}{}{\itshape}{}{\bfseries}{:}{.5em}{\thmnote{#3's }#1}
\theoremstyle{named}

\usepackage[Algorithm,ruled]{algorithm}
\usepackage{algpseudocode,algorithmicx}

\usepackage{cprotect} 
\usepackage{blkarray}

\newcommand{\RURVHaar}{\texttt{RURV\_Haar}}

\newcommand{\RURVROS}{\texttt{RURV\_ROS}}
\newcommand{\RVLUROS}{\texttt{RVLU\_ROS}}

\title{URV Factorization with Random Orthogonal System Mixing}
\author{Stephen Becker, James Folberth, Laura Grigori}

\begin{document}
\maketitle

\abstract{
The unpivoted and pivoted Householder QR factorizations are ubiquitous in numerical linear algebra.
A difficulty with pivoted Householder QR is the communication bottleneck introduced by pivoting.
In this paper we propose using random orthogonal systems to quickly mix together the columns of a matrix before computing an \emph{unpivoted} QR factorization.
This method computes a URV factorization which forgoes expensive pivoted QR steps in exchange for mixing in advance, followed by a cheaper, unpivoted QR factorization.
The mixing step typically reduces the variability of the column norms, and in certain experiments, allows us to compute an accurate factorization where a plain, unpivoted QR performs poorly.
We experiment with linear least-squares, rank-revealing factorizations, and the QLP approximation, and conclude that our randomized URV factorization behaves comparably to a similar randomized rank-revealing URV factorization, but at a fraction of the computational cost.
Our experiments provide evidence that our proposed factorization might be rank-revealing with high probability.
}

\section{Introduction}
\label{sec:intro}
The QR factorization of a matrix $A\in\reals^{m\times n}$ is a widely used decomposition, with applications in least-squares solutions to linear systems of equations, eigenvalue and singular value problems, and identification of an orthonormal basis of the range of $A$.
The form of the decomposition is $A=QR$, where $Q$ is $m\times m$ and orthogonal and $R$ is $m\times n$ and upper triangular.
When $A$ is dense and has no special structure, Householder reflections are often preferred to Gram-Schmidt (and its variants) and Givens rotations, due to their precise orthogonality and computational efficiency via the (compact) WY representation \cite{golub1998matrix, bischof1987wy, quintana1998blas}, which can utilize level-3 BLAS.
Indeed, Householder QR with a compact WY representation is implemented in the LAPACK routine \texttt{\_geqrf} \cite{anderson1999lapack}.

A common variant of the QR factorization is column pivoted QR, which computes the factorization $A\Pi = QR$, where $\Pi$ is a permutation matrix.
At the $i$th stage of the decomposition, the column of the submatrix $A(i:m,i:n)$ (in {\sc matlab} notation) with the largest norm is permuted to the leading position of $A(i:m,i:n)$ and then a standard QR step is taken.
The LAPACK routine \texttt{\_geqp3} implements column pivoted Householder QR using level-3 BLAS \cite{anderson1999lapack}.
However, it is typically much slower than the unpivoted \texttt{\_geqrf}, as \texttt{\_geqp3} still suffers from high communication costs \cite{demmel2015communication} and cannot be cast entirely in level-3 operations \cite{martinsson2015householder}.
We refer to Householder QR without pivoting as unpivoted QR (\texttt{QR}), and Householder QR with column pivoting as \texttt{QRCP}.

Improving on \texttt{QRCP}, recent works have used random projections to select blocks of pivots, emulating the behaviour of \texttt{QRCP}, while more fully utilizing level-3 BLAS \cite{duersch2015true, martinsson2015householder}.
Another approach uses so called ``tournament pivoting'' to select blocks of pivots and is shown to minimize communication up to polylogarithmic factors \cite{demmel2015communication}.
In each of these cases, a pivoted QR factorization is produced.

URV factorizations decompose $A$ as $A=URV$, where $U$ and $V$ have orthonormal columns and $R$ is upper triangular.
One can think of URV factorizations as a relaxation of the SVD, where instead of a diagonal singular value matrix, we require only that $R$ is upper-triangular.
Similarly, \texttt{QRCP} can be thought of as a URV factorization where $V$ is a permutation matrix, a special orthogonal matrix.
In Section \ref{sec:ls} we discuss how URV factorizations can be used to solve linear least-squares problems in much the same manner as QR factorizations or the SVD.

For example, let $V$ be a random orthogonal matrix sampled from the Haar distribution on orthogonal matrices.
The matrices $U$ and $R$ are computed with an \emph{unpivoted} QR factorization of $\hat{A}=AV^T$, and the resulting URV factorization is a strong rank-revealing factorization with high probability (see Subsection \ref{ssec:RURVHaar}) \cite{demmel2007fast}; we call this randomized factorization \RURVHaar.
This demonstrates that one can forego column pivoting at the cost of mixing together the columns of $A$ and still have a safe factorization.
However, taking $V$ to be a random, dense orthogonal matrix is not terribly computationally efficient, as $V$ is generated with an $n\times n$ unpivoted QR and must be applied with dense matrix multiplication.

We propose mixing with an alternating product of orthogonal Fourier-like matrices (e.g., discrete cosine, Hadamard, or Hartley transforms) and diagonal matrices with random $\pm 1$ entries, forming a so-called random orthogonal system (ROS) \cite{ailon2006approximate, tropp2011improved, mahoney2011randomized, meng2014lsrn}.
This provides mixing, but with a fast transform, as $V$ is never formed explicitly and can be applied with the FFT, or FFT-style algorithms (see Subsection \ref{ssec:RURVROS}).
We call this randomized URV factorization with ROS mixing \RURVROS.

\begin{figure}[ht!]
   \centering
   \includegraphics[width=0.8\textwidth]{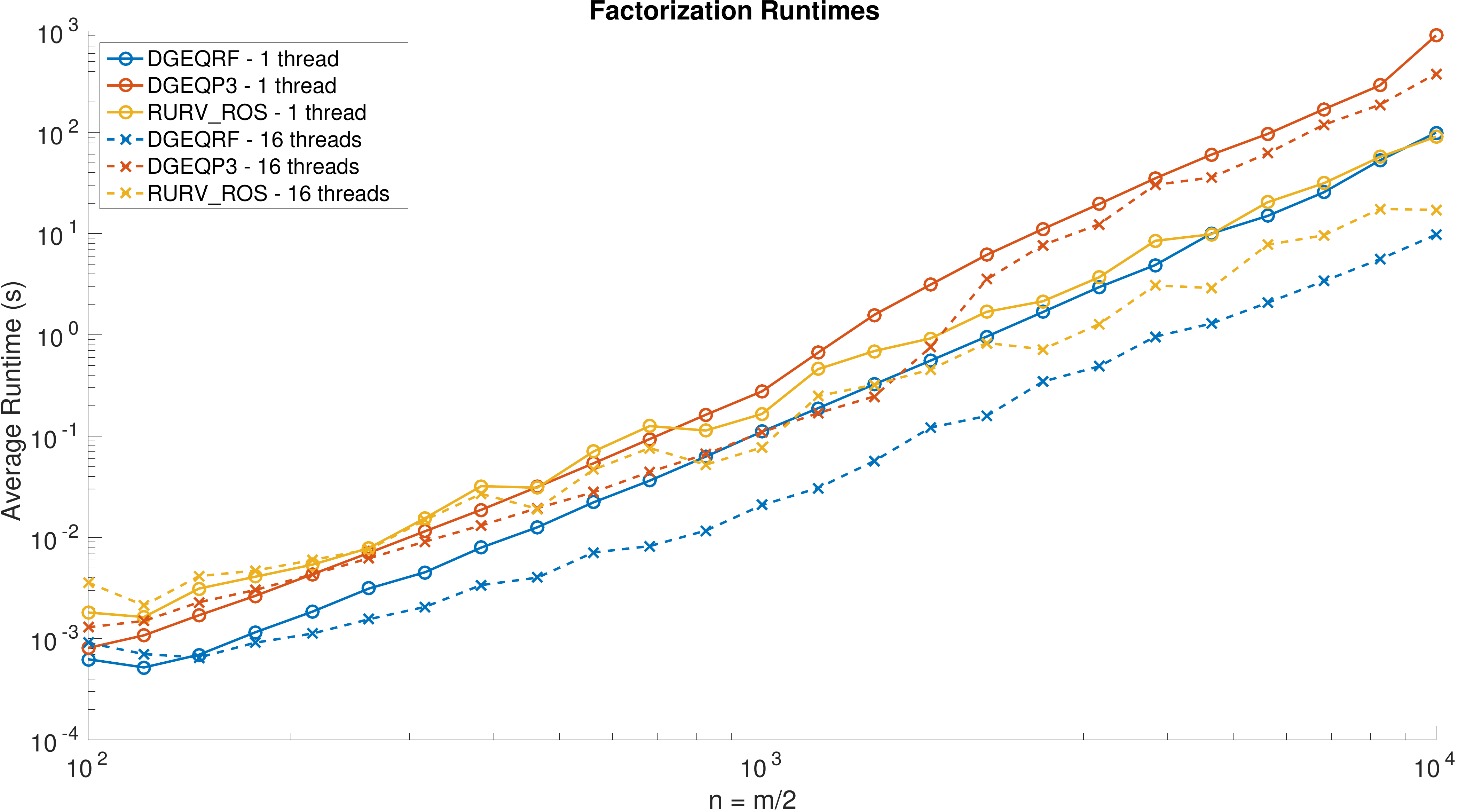}
   \caption{
      Average runtimes over five runs of \texttt{dgeqrf}, \texttt{dgeqp3}, and \RURVROS~on slightly tall-skinny matrices ($n=m/2$).
      Note that we do not include the time to generate the orthogonal factor $Q$ (labelled $U$ for \RURVROS), as all routines would use \texttt{dorgqr}.
      For the run with 16 threads, the sharp increase in runtimes beginning around size $2000\times1000$ matrices corresponds to the beginning of the regime where communication is the bottleneck of \texttt{QRCP}.
   }
   \label{fig:simple_qr_timing}
\end{figure}

Numerical experiments with our implementation of \RURVROS~demonstrate that for large matrices (i.e., where communication is the bottleneck of \texttt{QRCP}), \RURVROS~runs slightly slower than \texttt{\_geqrf} and significantly faster than \texttt{\_geqp3}.
Figure \ref{fig:simple_qr_timing} shows the average runtimes of \texttt{dgeqrf}, \texttt{dgeqp3}, and \RURVROS.
We used MATLAB's LAPACK \cite{matlab2016matlab} and the reference FFTW \cite{frigo2005design} with 1 and 16 threads on a desktop workstation with two Intel\textsuperscript{\textregistered} Xeon\textsuperscript{\textregistered} E5-2630 v3 CPUs running at $2.4$\,GHz.
See Subsection \ref{ssec:RURVROS} for more details on our implementation of \RURVROS.

Around $n=1000$, we begin to see a sharp increase in the runtime of \texttt{dgeqp3}, owing to the communication bottleneck of column pivoting.
In this region, \texttt{dgeqp3} with 16 threads does not see an appreciable improvement over running just a single thread.
In contrast, \texttt{dgeqrf} parallelizes much more nicely, as we can see an order of magnitude improvement in runtime when using 16 threads.
When using \RURVROS, we also see a noticeable improvement in runtime when using 16 threads versus 1 thread.

We also run timing and accuracy experiments on over- and underdetermined linear least-squares problems in Section \ref{sec:ls}.
In Subsection \ref{ssec:scaling_rr_cond} we sample the rank-revealing conditions of \cite{gu1996efficient,grigori2016low} for a variety of QR and URV factorizations, which suggest that \RURVROS~behaves similarly to \RURVHaar.
This provides evidence suggesting that \RURVROS~is rank-revealing with high probability.
We also examine using \RURVHaar~and \RURVROS~in a QLP approximation to the SVD in Subsection \ref{ssec:qlp_approx}.

\section{Randomized URV Factorization}
\subsection{Randomized URV Factorization via Haar Random Orthogonal Mixing}
\label{ssec:RURVHaar}

Demmel et al. proposed in \cite{demmel2007fast} a randomized URV factorization (RURV), which we call \RURVHaar, to use as part of eigenvalue and singular value decompositions.
Their RURV of an $m\times n$ matrix $A$ is based on sampling from the Haar distribution on the set of orthogonal (or unitary) matrices \cite{mezzadri2007generate}, using that sampled matrix to mix the columns of $A$, and then performing an \emph{unpivoted} QR on the mixed $A$, resulting in the factorization $A=URV$.

\begin{algorithm}[H]
   \caption{\texttt{RURV\_Haar} - Randomized URV with Haar mixing from \cite{demmel2007fast}}
   \begin{algorithmic}[1]
      \item[\textbf{Input:}] $A\in\reals^{m\times n}$
      \item[\textbf{Output:}] $U,R,V$
         \State Generate a random $n\times n$ matrix $B$ whose entries are i.i.d. $N(0,1)$.
         \State $[V,\hat{R}] = \texttt{qr}(B)$ \Comment{$V$ is Haar distributed; $\hat{R}$ is unused}
         \State $\hat{A} = AV^T$
         \State $[U,R] = \texttt{qr}(\hat{A})$
   \end{algorithmic}
   \label{alg:RURV_Haar}
\end{algorithm}

Haar orthogonal matrices are known to smooth the entries of the vectors on which they operate.
By multiplying $A$ on the right by a Haar orthogonal matrix $V^T$, we can mix together the columns of $A$, and reduce the variance of the column norms (see Figure \ref{fig:column_norms}).
The intuition behind the mixing is that by reducing the variance of the column norms, we reduce the effect that column pivoting would have, and can get away with unpivoted QR.
Indeed, in \cite{demmel2007fast} it is shown that Algorithm \ref{alg:RURV_Haar} produces a rank-revealing factorization with high probability, and can be used for eigenvalue and SVD problems.
It was further shown that Algorithm \ref{alg:RURV_Haar} produces a \emph{strong} rank-revealing factorization in \cite{ballard2010minimizing}.
Criteria for a (strong) rank-revealing factorization of the form $A=URV$ are as follows (taken from \cite{gu1996efficient,grigori2016low}, but slightly weaker conditions were used in \cite{ballard2010minimizing}):

\begin{enumerate}[1.]
   \item $U$ and $V$ are orthogonal and $R=\begin{bmatrix}R_{11}&R_{12}\\0&R_{22}\end{bmatrix}$ is upper-triangular, with $R_{11}$ $k\times k$ and $R_{22}$ $(n-k)\times (n-k)$;
      \item For any $1\le i \le k$ and $1\le j \le \min(m,n)-k$,
      
      \begin{equation}
         \label{eq:rr_cond_ratios}
          1\le \dfrac{\sigma_i(A)}{\sigma_i(R_{11})}, \dfrac{\sigma_j(R_{22})}{\sigma_{k+j}(A)}\le q(k,n),
      \end{equation}

         \noindent where $q(k,n)$ is a low-degree polynomial in $k$ and $n$.
      \item In addition, if 
		\begin{equation}
			\label{eq:rr_cond_strong}
			\|R_{11}^{-1}R_{12}\|_2
		\end{equation}
		\noindent is bounded by a low-degree polynomial in $n$, then the rank-revealing factorization is called \textbf{strong}.
\end{enumerate}

\noindent These conditions state that the singular values of $R_{11}$ and $R_{12}$ are not too far away from the respective singular values of $A$.
Thus, by performing a rank-revealing factorization instead of an expensive SVD, we can still gain insight into the singular values of $A$.

Both QR factorizations in Algorithm \ref{alg:RURV_Haar} are \emph{unpivoted}, and thus can be considerably cheaper than the standard column-pivoted Householder QR, \texttt{QRCP}.
However, a major drawback is the expense of generating and applying the random matrix $V$.
To sample an $n\times n$ matrix $V$ from the Haar distribution on orthogonal matrices, we take the $Q$ factor from an unpivoted QR factorization of an $n\times n$ matrix $B$ whose entries are i.i.d. $N(0,1)$ \cite{mezzadri2007generate}.
The dominant cost of this computation is the unpivoted QR factorization, which requires $\mathcal{O}(n^3)$ FLOPs.
We then compute $\hat{A}=AV^T$, which requires $\mathcal{O}(mn^2)$ FLOPs, followed by the unpivoted QR factorization to find $U$ and $R$, which costs $\mathcal{O}(mn^2)$ FLOPs.
To reduce the cost of forming and applying $V$, we propose replacing $V$ with a product of random orthogonal systems, which can each be applied implicitly and quickly, although providing slightly worse mixing.

\subsection{Randomized URV Factorization via Fast Random Orthogonal Mixing}
\label{ssec:RURVROS}

Consider a real $m\times n$ matrix $A$ and a product of random orthogonal systems (ROS) of the form

\begin{equation}
\label{eq:V_ROS}
V = \Pi\left[\prod_{i=1}^N FD_i\right],
\end{equation}

\noindent where each $D_i$ is a diagonal matrix of independent, uniformly random $\pm 1$ and $F$ is an orthogonal Fourier-like matrix with a fast transform.
Just like in \RURVHaar, we mix together the columns of $A$ as $\hat{A}=AV^T$.
The matrix $\Pi$ is a permutation matrix chosen so $\hat{A}\Pi^T$ sorts the columns of $\hat{A}$ in order of decreasing norm.
Replacing the Haar matrix $V$ in Algorithm \ref{alg:RURV_Haar} with the ROS based $V$ in \eqref{eq:V_ROS} yields the new algorithm we call \RURVROS, shown in Algorithm \ref{alg:RURV_ROS}.

\begin{algorithm}[H]
   \caption{\RURVROS~- Randomized URV with ROS mixing}
   \begin{algorithmic}[1]
   \item[\textbf{Input:}] $A\in\reals^{m\times n}$, number of mixing steps $N$, $\{D_i\}_{i=1}^N$ diagonal $\pm 1$ matrices
      \item[\textbf{Output:}] $U,R$ \Comment The $V$ matrix is not output because it is never explicitly formed
         \State $\hat{A} = A\prod_{i=N}^1(D_iF^T)$
         \State $\hat{A} = \hat{A}\Pi^T$  \Comment{Sort the columns of $\hat{A}$ so they are in order of decreasing $\ell_2$ norm.}
         \State $[U,R] = \texttt{qr}(\hat{A})$
         \State $V = \Pi\prod_{i=1}^N FD_i$
   \end{algorithmic}
   \label{alg:RURV_ROS}
\end{algorithm}

Each product $FD_i$ is referred to as a random orthogonal system (ROS) \cite{ailon2006approximate, tropp2011improved, mahoney2011randomized, meng2014lsrn}.
Examples of real-to-real, orthogonal Fourier-like transforms are the discrete cosine transform (e.g., DCT-II and DCT-III), the discrete Hartley transform, and the discrete Hadamard transform.
The Fourier-like matrix is never explicitly constructed, but rather is only used as an operator, for which we use a fast transform.
This brings the FLOP count for computing $AV^T$ from $\mathcal{O}(mn^2,n^3)$ to $\mathcal{O}(mn\log n)$.
In our experiments, we use the DCT-II and DCT-III for $F$ and $F^T$, as implemented in FFTW \cite{frigo2005design}.

\begin{figure}[ht]
   \centering
   \includegraphics[width=0.8\textwidth]{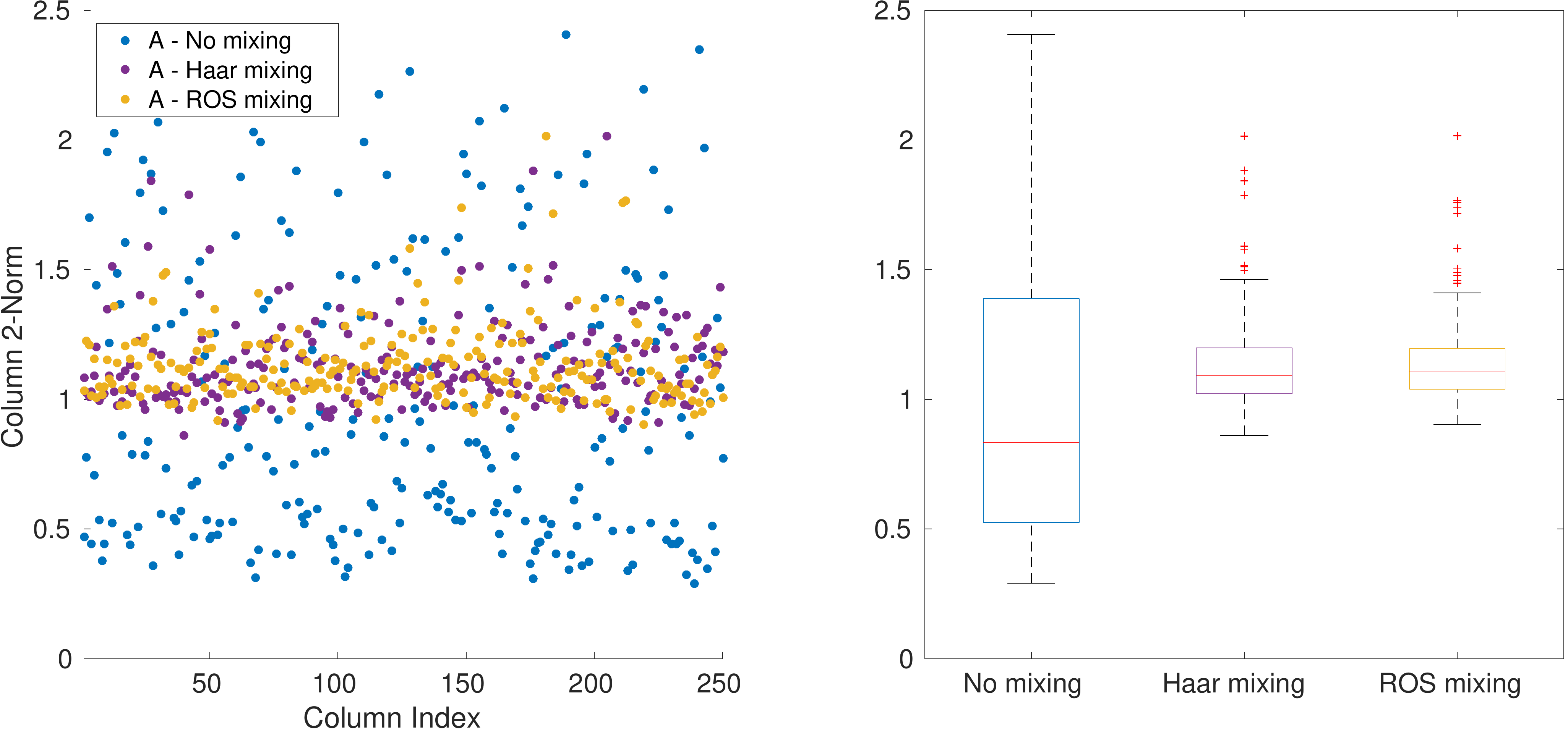}
   \caption{
      Mixing columns of $A$ together with Haar orthogonal matrices and ROS reduces the variance of column norms, while keeping the mean column norm about the same.
   }
   \label{fig:column_norms}
\end{figure}

Figure \ref{fig:column_norms} shows the effect of mixing with Haar matrices and ROS on the column norms of a random $250\times250$ matrix $A$, formed in {\sc matlab} with \verb|A = bsxfun(@times, randn(m,n)+exp(10*rand(m,n)), |\\\verb|exp(2*rand(1,n)))|, followed by \verb|A = A/mean(sqrt(sum(A.*A)))|, so that the mean column norm is one.
The variance of the column norms is clearly decreased by the mixing, and notably, Haar and ROS (with $N=1$) affect the distribution of column norms in a similar manner.
A theme of this paper is that \RURVROS~behaves similarly to \RURVHaar, which likely stems from their similar effect on the distribution of column norms.

To mix together the columns of $A$, we compute $\hat{A} = A\prod_{i=N}^1(D_iF^T)$.
The permutation/pre-sort matrix $\Pi$ is chosen so the columns of $\hat{A}\Pi^T$ are sorted in decreasing order of column norm.
The pre-sort is included to potentially enhance the accuracy and stability of \RURVROS.
The cost of this one-time, single sort is much smaller than the cost of the repeated column pivots in \texttt{QRCP}.

A {\sc matlab} implementation of \RURVROS~with $F$ taken to be the DCT-II is shown in Listing \ref{lst:rurv_ros.m}.
For in-core computations, it is sometimes more efficient to compute the mixing on left of $A^T$ via:

\[ AV^T = (VA^T)^T = \left(\Pi \prod_{i=1}^N \left(FD_i\right)A^T\right)^T. \]

\noindent This ``transpose trick'' is used in Listing \ref{lst:rurv_ros.m} for efficiency, and also to cleanly interface with {\sc matlab}'s \texttt{dct} function, which applies the transform to the columns of its input.
Listing \ref{lst:rurv_ros.m} explicitly returns $U$ and $R$ from the factorization, but returns function handles for $V$ and $V^T$, which can be used to apply $V$ and $V^T$, respectively, to the left side of their input.

The implementation used for our experiment is similar, but has performance-critical sections written in C using {\sc matlab}'s MEX interface.
The mixing step is performed in C using FFTW and the unpivoted QR is performed in C using LAPACK routines from {\sc matlab}'s LAPACK \cite{frigo2005design,anderson1999lapack, matlab2016matlab}.
The use of FFTW gives us great control over how the transform is applied (e.g., in blocks, multithreaded, perhaps not utilizing the ``transpose trick'', etc.).
More details on the use of FFTW for mixing are given in Subsection \ref{ssec:ls_timing}.

\begin{lstlisting}[language=Matlab, caption=A {\sc matlab} implementation of \RURVROS, label=lst:rurv_ros.m]
function [U,R,V,Vt] = RURV_ROS(A, n_its)
% RURV_ROS  RURV with ROS mixing for real matrices

   [m,n] = size(A);
   D_diags = sign(rand(n,n_its)-0.5); % diagonals of D_i; i.i.d. uniform +- 1
   
   % ROS mixing Ahat = A*V' using "transpose trick"
   Ahat = apply_V(A',D_diags)';

   % pre-sort
   nrms = sqrt(sum(Ahat.^2,1));
   [nrms,p] = sort(nrms, 2, 'descend');
   p_inv(p) = 1:numel(p);
   Ahat = Ahat(:,p);

   % unpivoted QR factorization
   [U,R] = qr(Ahat,0);
   
   % Return function handles to apply V and V^T on the left
   V = @(A) apply_V(A,D_diags,p);
   Vt = @(A) apply_Vt(A,D_diags,p_inv);
end

function [Ahat] = apply_V(A, D_diags, p)
% apply_V  Apply ROS mixing: V*A
   Ahat = A;
   for i=1:size(D_diags,2)
      Ahat = bsxfun(@times, Ahat, D_diags(:,i)); % Ahat = D_i*Ahat
      Ahat = dct(Ahat);                       % Ahat = F*Ahat;
   end
   if nargin == 3 % apply sorting
      Ahat = Ahat(p,:);
   end
end

function [Ahat] = apply_Vt(A, D_diags, p_inv)
% apply_Vt  Apply transpose ROS mixing: V^T*A
   if nargin == 3 % apply sorting
      Ahat = A(p_inv,:);
   else
      Ahat = A;
   end
   for i=size(D_diags,2):-1:1
      Ahat = idct(Ahat);                      % Ahat = F^T*Ahat;
      Ahat = bsxfun(@times, Ahat, D_diags(:,i)); % Ahat = D_i*Ahat
   end
end
\end{lstlisting}

\section{Applications to Least-Squares Problems}
\label{sec:ls}

\subsection{Solving Least-Squares Problems with a URV Factorization}

A URV factorization can be used to solve least-squares problems in much the same manner as a QR factorization.
Throughout this subsection we assume that $A$ is $m\times n$ and full-rank.
We are interested in finding a solution to

\[ \min_x \|Ax-b\|_2 \] 

\noindent for both the overdetermined case $m\ge n$ and the underdetermined case $m < n$.

\subsubsection{Overdetermined Systems}
Consider first the case when $A$ is overdetermined.
To find the least-squares solution with a QR factorization, we only need a thin QR factorization, where $Q$ is $m\times n$ and $R$ is $n\times n$ \cite{golub1998matrix}.
Similarly, the internal QR factorization in \RURVROS~can be a thin QR.
By computing $A=URV$ and using that $U$ has orthonormal columns,

\[ \min_x \|Ax-b\|_2 = \min_x \|URVx-b\|_2 = \min_x \|RVx-U^Tb\|. \] 

\noindent The least-squares problem reduces to the non-singular $n\times n$ upper-triangular system $Ry=U^Tb$ in the auxiliary variable $y=Vx$.
The system $Ry=U^Tb$ is solved implicitly for $y$ with backward substitution, and then the least-squares solution is found with $x=V^Ty$.

Note that we do not need to explicitly form $U$ to apply $U^T$ to $b$.
When we call LAPACK's \texttt{\_geqrf} on $\hat{A}=AV^T$, the routine overwrites the upper-triangular part of $\hat{A}$ with $R$ and the Householder reflectors in the strictly lower triangular part of $\hat{A}$.
By feeding the Householder reflectors into \texttt{\_ormqr}, we can implicitly compute $U^Tb$ in $\mathcal{O}(mn)$ FLOPs without ever accumulating $U$ \cite{anderson1999lapack}.

The dominant cost of using \RURVHaar~to compute least-squares solutions is a mix of generating $V$, computing $\hat{A}=AV^T$, and the thin QR to find $U$ and $R$.
The latter two operations cost $\mathcal{O}(mn^2)$ FLOPs.
The dominant cost of using \RURVROS~is also $\mathcal{O}(mn^2)$, but the leading cost term only comes from the unpivoted QR, as the mixing $\hat{A}=AV^T$ is $\mathcal{O}(mn\log n)$.

\subsubsection{Underdetermined Systems}
Now consider the underdetermined case.
A full URV factorization $A=URV$ is of the following form:

\begin{equation}
   \label{eq:underdetURV}
   \begin{blockarray}{cc}
      & \qquad n\qquad\\
      \begin{block}{c[c]}
         m & A\\
      \end{block}
   \end{blockarray}
   =
   \begin{blockarray}{cc}
      & m\\
      \begin{block}{cc}
         m & U\\
      \end{block}
   \end{blockarray}
   \begin{blockarray}{cc}
      m & n-m\\
      \begin{block}{[cc]}
         R_{11} & R_{12}\\
      \end{block}
   \end{blockarray}
   \begin{blockarray}{cc}
      n&\\
      \begin{block}{cc}
         V&n\\
      \end{block}
   \end{blockarray}
\end{equation} 

\noindent Since $A$ is assumed to be full-rank, $\min_x \|Ax-b\|_2=0$ and we seek to solve $Ax=b$.
As in the overdetermined case, make the change of variable $y=Vx$; we now consider solving the upper-trapezoidal system $Ry=U^Tb$.
Partitioning $y$ into $m\times 1$ and $(n-m)\times 1$ blocks results in the block system

\[ \begin{bmatrix} R_{11} & R_{12}\end{bmatrix}\begin{bmatrix} y_1\\y_2\end{bmatrix} = U^Tb, \] 

\noindent where $R_{11}$ is upper-triangular and full-rank.
A particularly simple solution is found by setting $y_2=0$ and performing backward substitution to find $y_1$.
Following \cite{golub1998matrix}, we call this the \textbf{basic solution}.
Note that the basic solution has $m-n$ zeros in $y$, but after unmixing to find $x_\text{basic}=V^Ty$, the zeros in $y_2$ are mixed with the nonzeros in $y_1$, destroying the sparsity of $x_\text{basic}$.
While this is less than ideal, mixing and unmixing is fast, and sparsity in the mixed domain might still be applicable in certain problems.

Notice that $R_{12}$ is not used to compute the basic solution.
Since $R$ is computed from $\hat{A}=AV^T$, which mixes all the columns of $A$ together, we may compute $U$ and $R_{11}$ from the $QR$ factorization of $\hat{A}(:,1:m)$ (in {\sc matlab} notation).
This avoids the computation of $R_{12}$, leading to a faster solution.
Mixing to find $\hat{A}=AV^T$ costs $\mathcal{O}(mn\log n)$; computing $R_{11}$ costs $\mathcal{O}(m^3)$; and applying $U^Tb$, backward substitution to find $y=R_{11}^{-1}U^Tb$, and unmixing to find $x_\text{basic}=V^Ty$ all cost a negligible amount for large $m$ and $n$.
This brings the total cost to compute the basic solution to $\mathcal{O}(m^3, mn\log n)$ FLOPs.

Another common solution is the \textbf{minimum norm solution}.
Since the solution set $\mathcal{X}=\{x\in\reals^n \,|\, Ax=b\}$ is closed and convex, there exists a unique minimum norm solution, which is a principal attraction to the minimum norm solution (a similar statement holds even when $A$ is rank deficient).
Finding the minimum norm solution can be expressed as the problem

   \[ \begin{array}{ll}\min & \|x\|^2\\\text{s.t.} & Ax=b.\end{array} \] 

\noindent Let $\mathcal{L}(x,\nu) = x^Tx + \nu^T(Ax-b)$ be the Lagrangian function.
Slater's condition for this problem is simply that the problem is feasible, which is of course satisfied since we assume $A$ is full-rank.
Therefore, strong duality holds and the KKT conditions,
   
  \[ \nabla_x \mathcal{L} = 2x + A^T\nu = 0, \quad Ax-b=0, \]

\noindent give necessary and sufficient conditions for the solution \cite{boyd2004convex}.
Solving the KKT conditions gives $x_\text{mn}=A^T(AA^T)^{-1}b = A^\dagger b$, where $A^\dagger$ is the (right) pseudoinverse of $A$.
To use this closed-form solution efficiently, it is convenient to perform a QR factorization of $A^T$.
Specifically, if we let $A^T=QR$, then $x_\text{mn} = QR^{-T}b$, where $R^{-T}b$ is computed implicitly with forward substitution.

To find the minimum norm solution with mixing, we should mix the columns of $A^T$ in preparation for the unpivoted QR of $\hat{A}^T$.
Let $\hat{A}^T=A^TV^T$ (which we may compute via $\hat{A}=VA$) and compute $\hat{A}^T=U^TL^T$ via unpivoted QR.
We then have the factorization $A=V^TLU$, where $V$ is our fast ROS mixing matrix, $L$ is $m\times m$ lower triangular, and $U$ is $m\times n$ with orthonormal rows (i.e., $U^T$ is orthonormal).
We call the algorithm to compute $A=V^TLU$ \RVLUROS~in analogy with \RURVROS.
By multiplying $Ax=b$ on the left by $V$, we find $\hat{A}x=Vb$, and from the discussion above, the minimum norm solution is $x_\text{mn}=U^TL^{-1}Vb$.
Again note that $L^{-1}$ is applied implicitly using forward substitution.
The dominant cost of this approach is again the unpivoted QR factorization of $\hat{A}^T$, which costs $\mathcal{O}(mn^2)$ FLOPs, which can be significantly higher than the $\mathcal{O}(m^3,mn\log n)$ FLOPs for the basic solution. 

\subsection{Timing Experiments}
\label{ssec:ls_timing}

Solving least-squares problems with \RURVROS~or \RVLUROS~factorizations will be slightly slower than using unpivoted QR; the additional cost comes almost entirely from the mixing steps in Algorithm \ref{alg:RURV_ROS}.
In our code, we use the DCT-II and DCT-III, as implemented in FFTW \cite{frigo2005design}.
For improved performance, we cache FFTW ``wisdom'' in a file and load it the next time it is applicable.
Finding the solution proceeds in three stages: mixing to find $\hat{A}$, performing unpivoted QR factorization of $\hat{A}$ or $\hat{A}^T$, and computing the final solution vector, which may involve mixing a single vector.
For moderately large overdetermined problems, mixing to find $\hat{A}$ takes about $25\%$ of the total runtime; unpivoted QR factorization $75\%$ of the total time; and solving/mixing takes a negligible amount of time, since it is applied to only a single vector.

We compare with BLENDENPIK, which uses mixing across rows and row sampling to find a good preconditioner for LSQR \cite{avron2010blendenpik, paige1982lsqr}.
The authors wrote most of their code in C for efficiency, calling LAPACK and FFTW libraries and providing their own implementation of LSQR.
When we installed BLENDENPIK, we precomputed FFTW ``wisdom'' for the most patient planner setting, which results in the highest performance at run-time.
In the underdetermined case, BLENDENPIK computes the minimum norm solution.
With the exception of using DCT mixing, we used the default parameters provided in BLENDENPIK's interface.

It is worth noting that the well-known backslash (\texttt{\textbackslash}) operator in {\sc matlab} solves (rectangular) linear systems in the least-squares sense using a QR-based approach.
{\sc matlab}'s \texttt{\textbackslash} operator tends to be significantly slower than BLENDENPIK and \RURVROS, but \texttt{\textbackslash} also supports the case of rank-deficient matrices \cite{matlab2016matlab}.
LAPACK has a variety of least-squares routines, and can handle full-rank and rank-deficient matrices.
The LAPACK routine \texttt{\_gels} uses a simple rescaling and unpivoted QR or LQ to solve full-rank least-squares problems \cite{anderson1999lapack}.
For highly overdetermined systems, BLENDENPIK is reported to beat QR-based solvers, including \texttt{\_gels}, by large factors \cite{avron2010blendenpik}.

For the following timing experiments, we take $A$ to be a random matrix constructed by $A=U\Sigma V^T$ where $U$ and $V$ are random orthogonal matrices and $\Sigma$ is a diagonal matrix of singular values such that $\kappa_2(A)=10^6$ ($\kappa_2(A)=\|A\|_2\|A^{-1}\|_2$ is the spectral condition number of $A$).
We take a single random right-hand side vector $b$ with entries sampled from $N(0,1)$ and solve the problem $\min_x \|Ax-b\|_2$.
We link BLENDENPIK and our code against {\sc matlab}'s LAPACK and the standard FFTW library.
For timing results, we run each routine once to ``warm-up'' any JIT-compiled {\sc matlab} code, and run a number of samples.  

Our code is designed to scale up to multiple threads on a single machine, using multi-threaded versions of LAPACK and FFTW, but BLENDENPIK currently uses only a single thread for their FFTW calls.
We note that it would be straightforward to extend BLENDENPIK to use multi-threaded FFTW calls, but mixing is hardly the dominant cost of BLENDENPIK, so one would not expect see a large improvement in runtimes.
Nevertheless, we perform the following timing experiments using only a single thread in order to compare fairly BLENDENPIK and \RURVROS.

\begin{figure}[h!]
   \centering
   \includegraphics[width=0.49\textwidth]{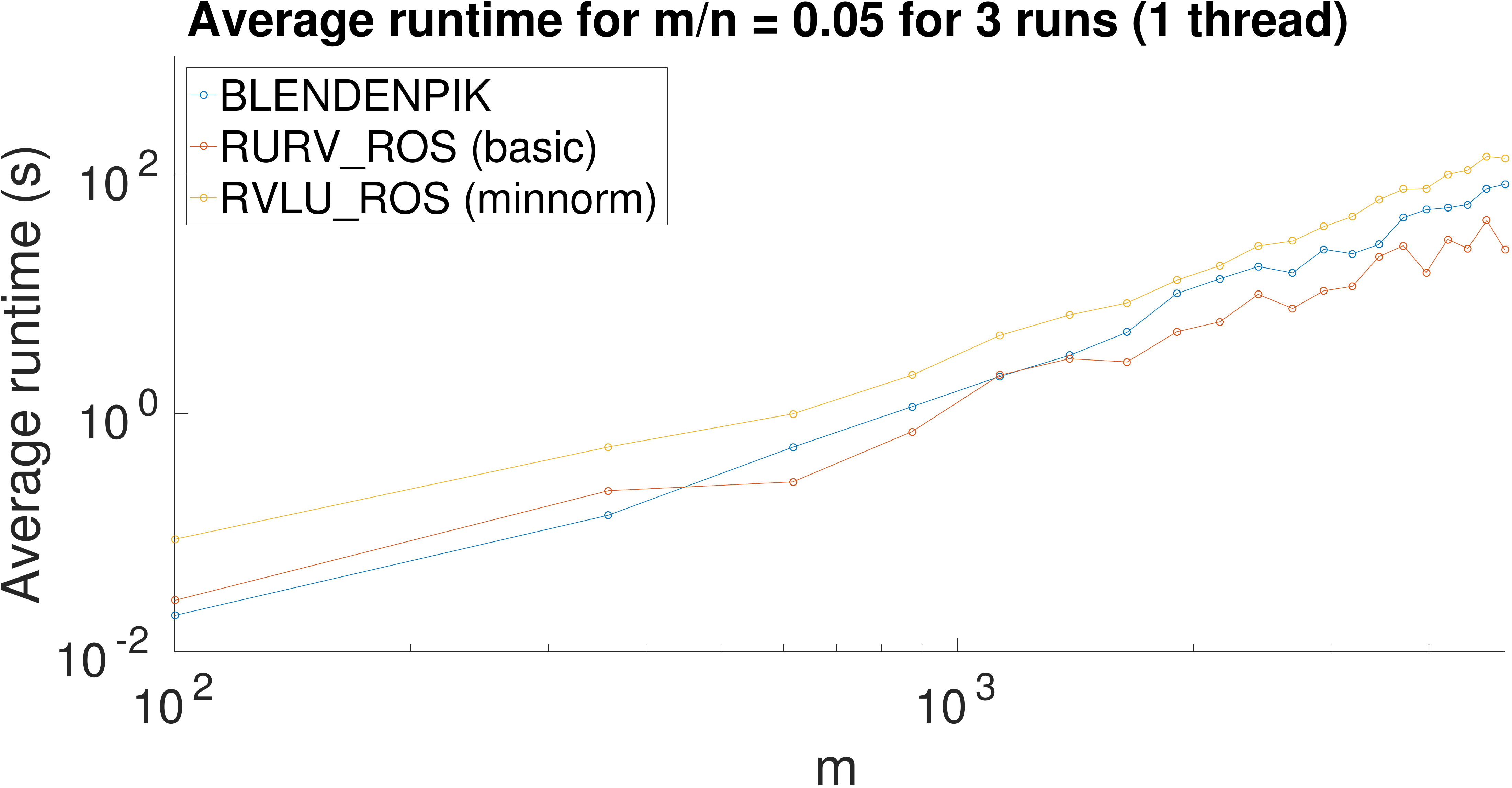}
   \includegraphics[width=0.49\textwidth]{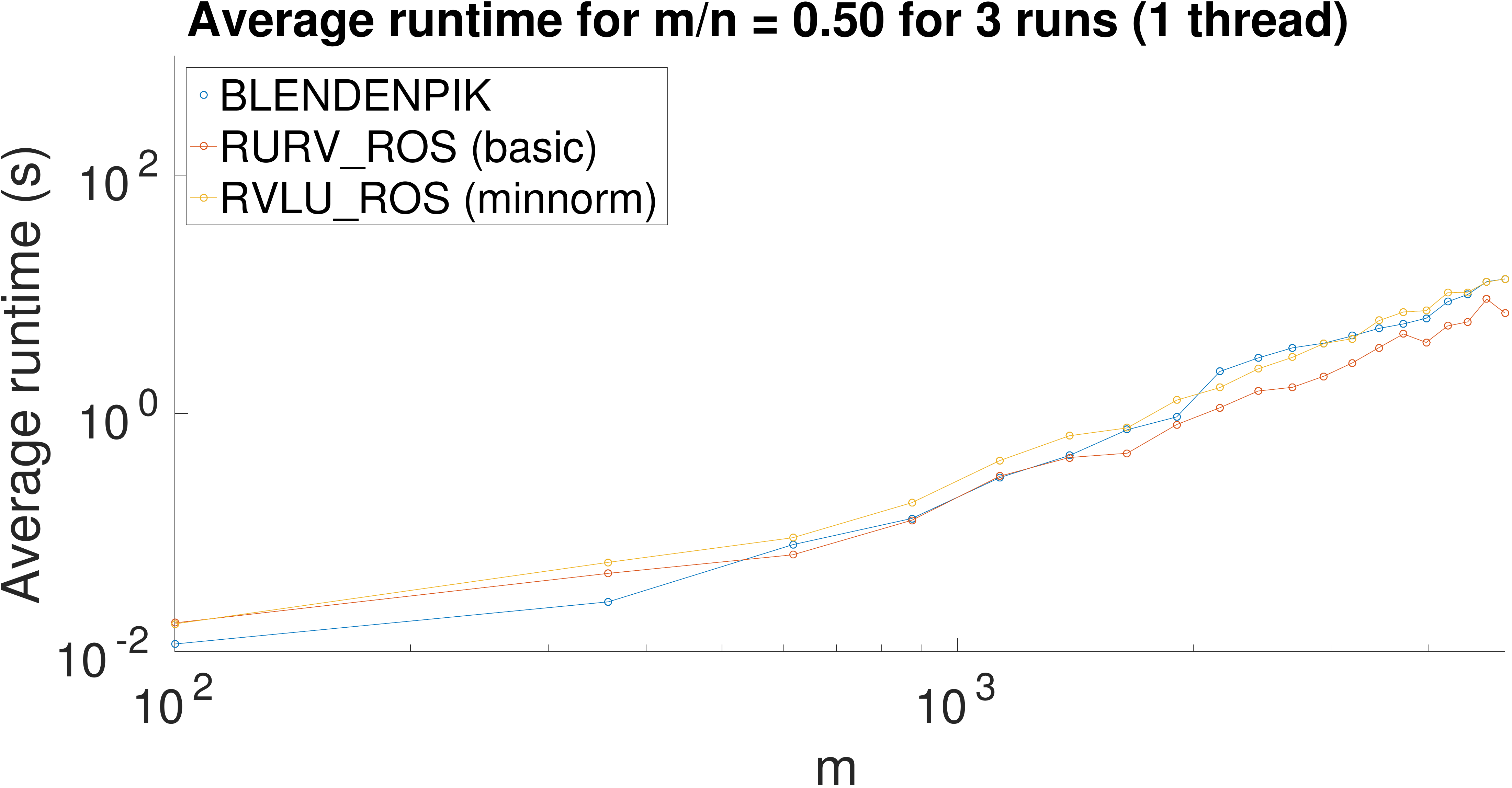}
   \caption{
      Average runtime for BLENDENPIK, \RURVROS, and \RVLUROS~approaches on moderately and slightly underdetermined systems.
      The \RVLUROS~based minimum norm solution is consistently slower than BLENDENPIK, which also computes the minimum norm solution.
      The basic solution computed with \RURVROS, being simpler to compute, is a considerably faster than the \RVLUROS~based minimum norm solution.
   }
   \label{fig:times_underdet}
\end{figure}

\begin{figure}[h!]
   \centering
   \includegraphics[width=0.49\textwidth]{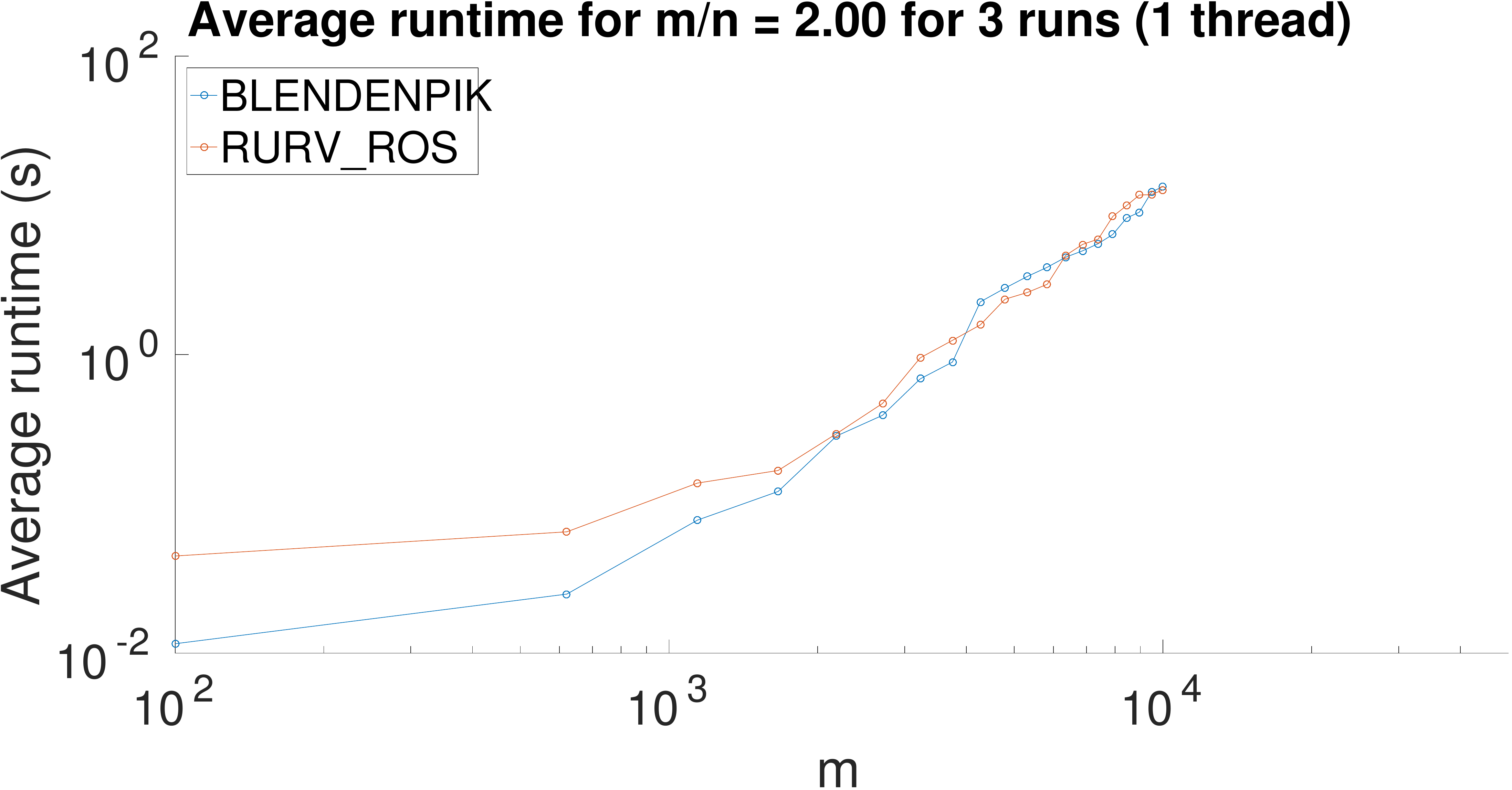}
   \includegraphics[width=0.49\textwidth]{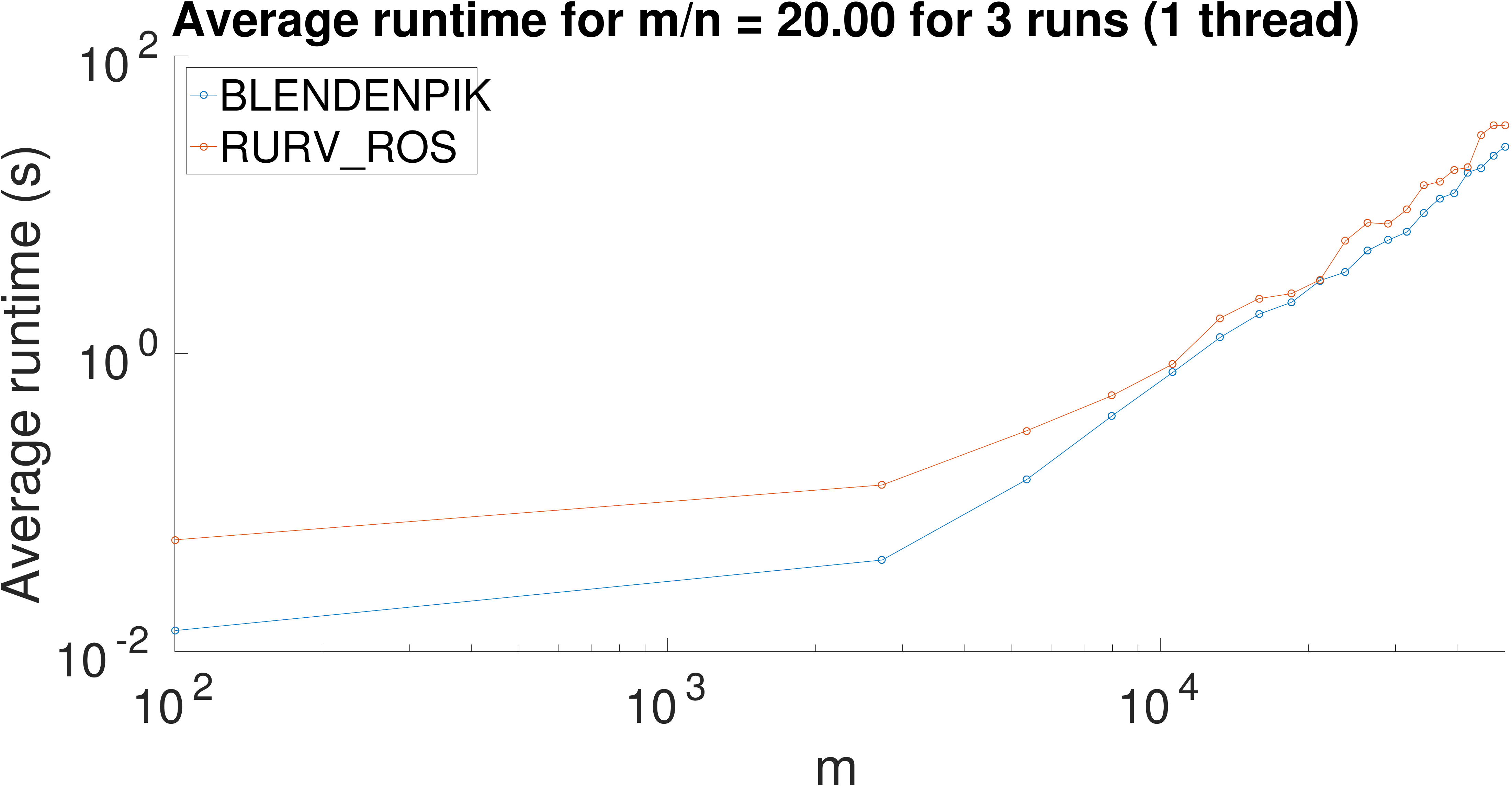}
   \caption{
      Average runtime for BLENDENPIK and the \RURVROS~approach on overdetermined systems.
      \RURVROS~compares favorably to BLENDENPIK for $m/n$ close to one.
      For larger $m/n$ (i.e., more highly overdetermined), BLENDENPIK performs better than \RURVROS, which we expect from \cite{avron2010blendenpik}, which shows that BLENDENPIK tends to outperform QR-based solvers on highly overdetermined systems.
   }
   \label{fig:times_overdet}
\end{figure}

Figures \ref{fig:times_underdet} and \ref{fig:times_overdet} show the average runtime for random $A$ with $\kappa_2(A)=10^6$ of various sizes.
We consider moderately underdetermined, slightly underdetermined, slightly overdetermined, and moderately overdetermined examples.
For underdetermined systems, computing the basic solution with \RURVROS~is slightly faster than BLENDENPIK, which computes the minimum norm solution.
Using \RVLUROS~to compute the minimum norm solution is moderately slower than BLENDENPIK.

Notice in Figure \ref{fig:times_underdet} that the average runtime for the basic solution is slightly jagged.
This variance is due to FFTW's planner finding a faster plan for certain sizes.
We could potentially improve the runtime by optionally zero-padding $A$ and using transforms of a slightly larger size, which may allow FFTW's planner to find a faster plan for the larger size.
Using zero-padding would change the values in the mixed matrix $\hat{A}$, however, so for now we do not investigate using zero-padding.

\begin{figure}[h!]
   \centering
   \includegraphics[width=0.49\textwidth]{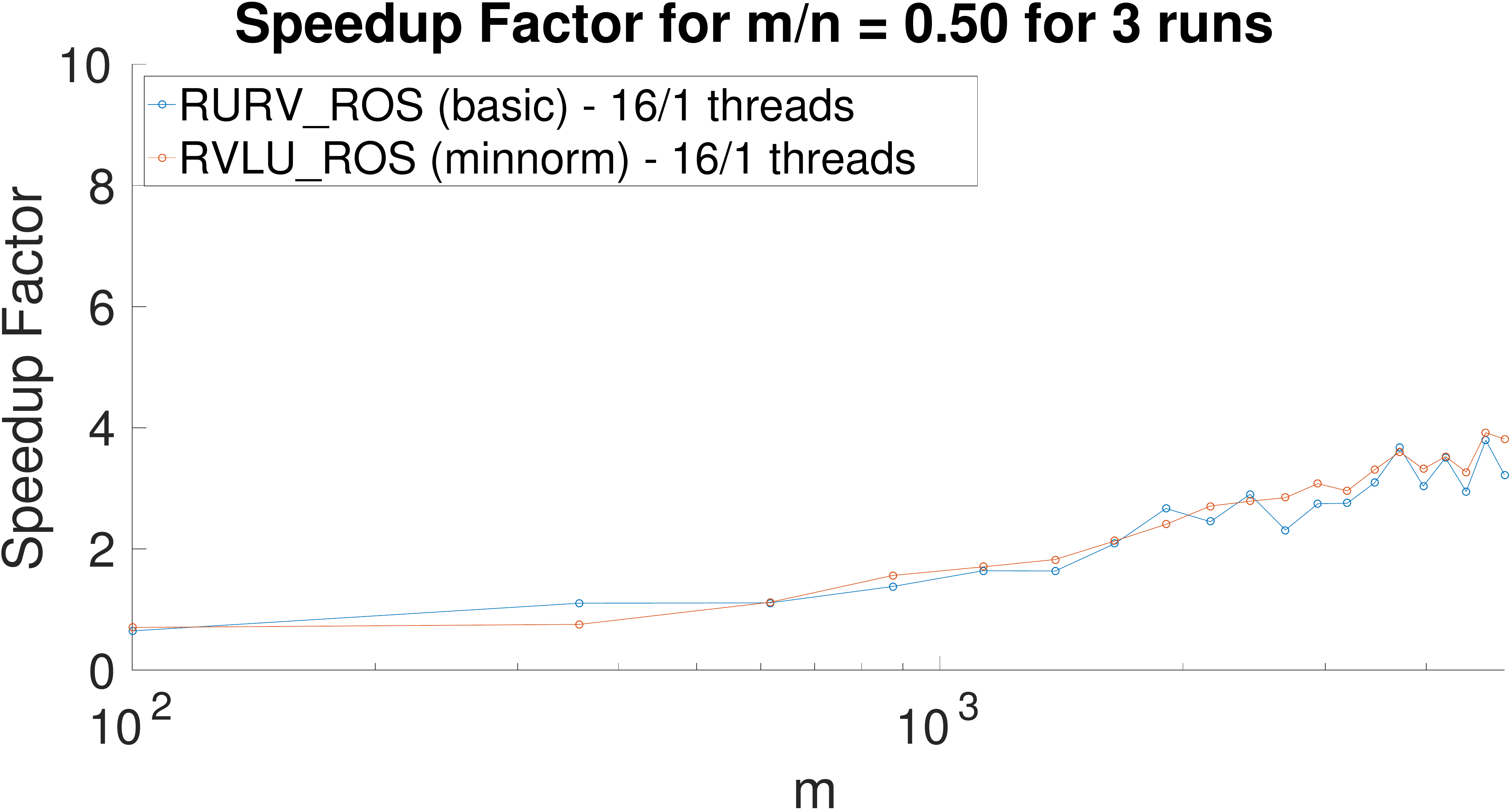}
   \includegraphics[width=0.49\textwidth]{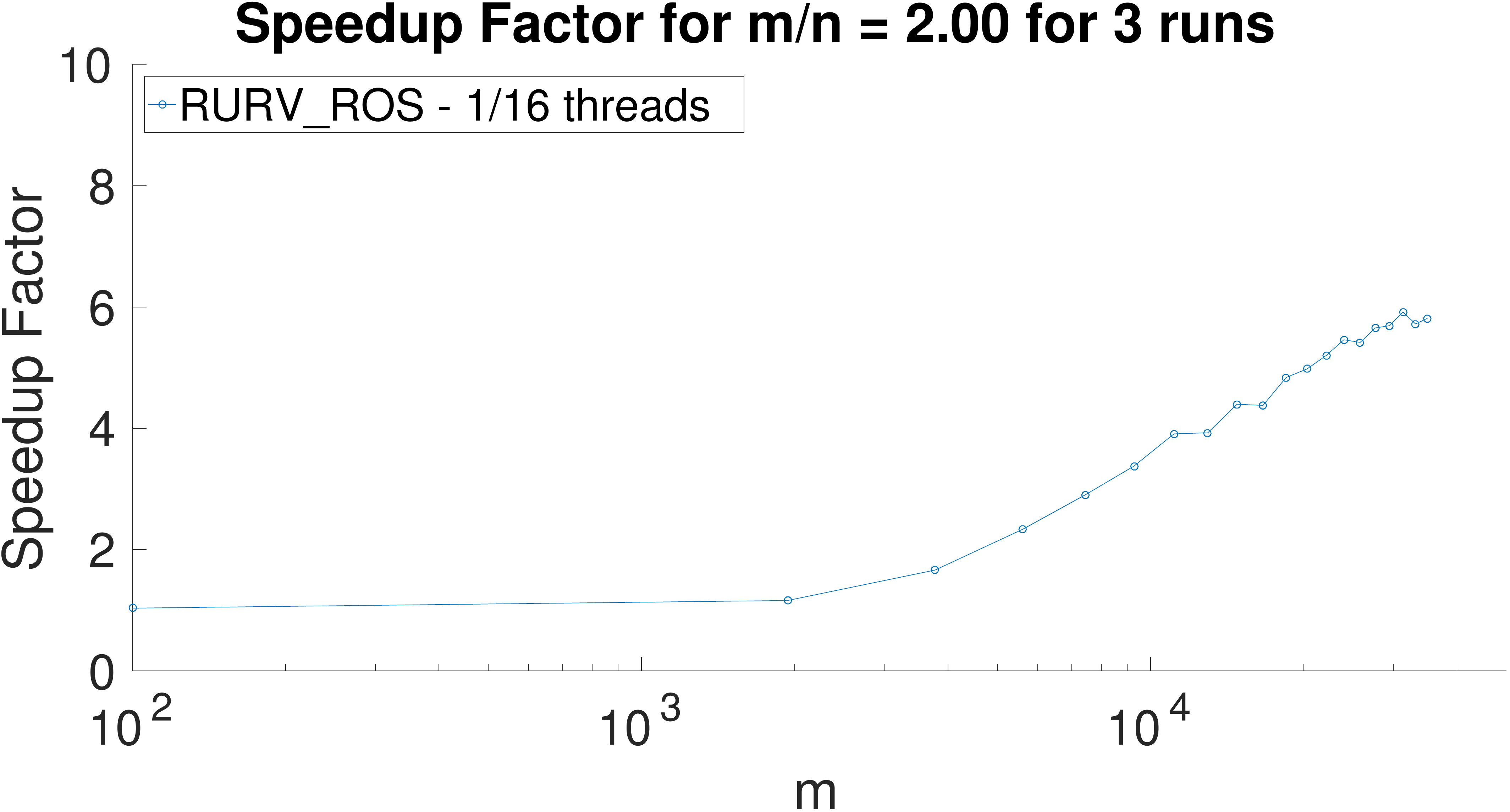}
   \caption{
      Ratios of average runtime for the \RURVROS~and \RVLUROS~approaches on slightly under- and overdetermined systems using 1 and 16 threads.
      We plot the ratio of the runtime using 1 thread over the runtime using 16 threads, so speedup factors greater than 1 correspond to an improvement when running in parallel.
   }
   \label{fig:times_pair}
\end{figure}

Figure \ref{fig:times_pair} shows the ratio of the runtimes for \RURVROS~and \RVLUROS~for both 1 and 16 threads.
For slightly underdetermined systems, the speedup factor approaches 4 (i.e., for large $m$, using 16 threads runs about 4 times faster than using only 1 thread).
For slightly overdetermined systems, the speedup factor increases for $m\ge 1000$ and approaches 6.
Although we see speedup factors of less than ideal 16, our implementation does parallelize nicely.
The speedup factors may very well continue to increase outside of the range of matrices we tested (until we run out of memory on the machine, that is).

The larger speedup factor for overdetermined systems is likely due to the ratio of the work computed in the mixing stage and the factorization stage.
For overdetermined systems, mixing occurs along the smaller dimension of the matrix, so there are many smaller transforms, compared to underdetermined systems.
This gives a higher proportion of the work to the QR factorization, during which we can more effectively utilize additional cores.

Before we continue with our discussion of our experiments, we make a brief note on implementing \RURVROS~on a distributed memory machine.
The two major steps of \RURVROS~are mixing and the unpivoted QR factorization, which can be handled by FFTW and ScaLAPACK, respectively.
FFTW has a distributed memory implementation using MPI, and the interface is very similar to the shared memory interface.
ScaLAPACK's routine \texttt{p\_geqrf} performs unpivoted Householder QR, but uses a suboptimal amount of communication.
In \cite{demmel2012communication} communication-avoiding QR (CAQR) was introduced.
CAQR sends a factor of $\mathcal{O}(\sqrt{mn/P})$ fewer messages than \texttt{p\_geqrf} (where $P$ is the total number of processors in the grid), and achieves the optimal amount of communication (up to polylogarithmic factors).
The reduction in communication is predicted to result in significantly faster factorization runtimes.
Using FFTW, ScaLAPACK, or other, existing codes as building blocks, we expect that \RURVROS~can be implemented efficiently and straightforwardly for distributed memory environments.

\subsection{Example - Correlated Columns and the Basic Solution}
The dominant cost of using a URV factorization to compute the basic solution to an underdetermined $m\times n$ system ($m < n$) is computing a QR factorization of $\hat{A}(:,1:m)$.
Thus, it is asymptotically cheaper than the minimum norm solution, which uses an LQ factorization of the full $\hat{A}$.
Specifically, computing the basic solution costs $\mathcal{O}(m^3,mn\log n)$ FLOPs, while computing the minimum norm solution costs $\mathcal{O}(m^2n)$ FLOPs \cite{golub1998matrix}.
We have seen in Figure \ref{fig:times_underdet} that computing the basic solution is significantly faster than computing the minimum norm solution with \RURVROS, and even slightly faster than BLENDENPIK, which also computes the minimum norm solution.

The following simple example shows that finding the basic solution using unpivoted QR is numerically unstable for some least-squares problems.
Consider

\[ A = \begin{bmatrix}1&0&0&0\\0&1&1&1\\0&0&\epsilon&1\end{bmatrix}, \] 

\noindent with $\epsilon \ll 1$ and consider solving $\min_x \|Ax-b\|_2$.
Note that $A$ is full-rank and that it can be analytically verified that $\kappa_2(A) \sim 1+\sqrt{2}$ as $\epsilon\to 0$, so $A$ is very well conditioned for $\epsilon\ll 1$.
However, columns 2 and 3 are increasingly correlated as $\epsilon\to 0$.
Since $A$ is already upper trapezoidal, an unpivoted QR factorization does not change $A$ when finding $R$.
When finding the basic least-squares solution in this manner, it transpires that we solve the linear system

\[ R_{11}x_1 = \begin{bmatrix}1&0&0\\0&1&1\\0&0&\epsilon\end{bmatrix}\begin{bmatrix}x_{1,1}\\x_{1,2}\\x_{1,3}\end{bmatrix} = Q^Tb = b = \begin{bmatrix}b_1\\b_2\\b_3\end{bmatrix}. \] 

\noindent However, this system has $\kappa_2(R_{11})\sim 2/\epsilon$ as $\epsilon\to 0$, and so can become quite ill-conditioned as columns 2 and 3 become more correlated.
Even for this small system, finding the basic solution $\hat{x}$ with an unpivoted QR factorization leads to a large residual $\|A\hat{x}-b\|_2$.
   
This can be fixed by using \texttt{QRCP} on all of $A$ instead of unpivoted QR on $A(:,1:m)$.
Note that using \texttt{QRCP} on $A(:,1:m)$ will encounter similar ill-conditioning problems, as doing so will not allow \texttt{QRCP} to pivot in column 4.
Using a URV factorization which mixes in column 4 of $A$ also leads to a much better conditioned $R_{11}$, alleviating the issue of correlated columns when using the basic solution.
This simple example can be extended to larger matrices, as we show next.

Consider an $m\times (n-p)$ matrix $A$ with $m < n$ such that each element of $A$ is sampled from $N(0,1)$.
We augment $A$ by adding $p$ randomly selected columns of $A$ to the end of $A$, making $A$ $m\times n$.
The augmented $A$ has $p$ perfectly correlated columns, so we add a small amount of $N(0,\sigma^2)$ noise to the augmented $A$ so the correlation is not perfect.
We then randomly shuffle the columns.
Listing \ref{lst:corr_mat} gives {\sc matlab} code to generate such a matrix, which tends to be well-conditioned.
If the final permutation places a pair of highly correlated columns in the first $m$ columns of $A$, finding the basic solution $\hat{x}$ with unpivoted QR will produce an ill-conditioned $R_{11}$, leading to a large residual $\|A\hat{x}-b\|_2$.
This can be solved by mixing with \RURVHaar~or \RURVROS, or computing the (more costly) minimum norm solution.

\begin{lstlisting}[language=Matlab,caption={\sc matlab} code to generate a matrix with a few correlated columns,label=lst:corr_mat]
m = 1000; n = 1500;
p = 10; e = 1e-4;

A = randn(m,n-p);
perm = randperm(n-p,p);       
A = [A A(:,perm)];
perm = randperm(n);           
A = A(:,perm) + e*randn(m,n); 
\end{lstlisting}

Table \ref{tab:basic_resids} shows the residuals and runtimes on a matrix generated with Listing \ref{lst:corr_mat}.
We tested unpivoted \texttt{QR}, \texttt{QRCP}, BLENDENPIK, \RURVHaar, and \RURVROS.
We use unpivoted \texttt{QR} on only the first $m$ columns of $A$, so it produce a significantly larger residual than the other methods.
Note that BLENDENPIK actually computes the minimum norm solution and is included for reference.
\RURVHaar~and \RURVROS~both compute basic solutions with acceptably small residuals.
As expected, \RURVROS~is considerably faster than \RURVHaar, but slightly slower than plain, unpivoted \texttt{QR}.

It is interesting to note that the norm of the mixed basic solution is considerably smaller than the unmixed basic solution.
Table \ref{tab:basic_resids2} shows the same comparison for $p=0$ correlated columns, where we see that the mixed and unmixed basic solutions have norms that are not unreasonably large.
The norms of the mixed basic solutions are on the same order for the cases of $p=10$ and $p=0$ correlated columns, unlike \texttt{QR} and \texttt{QRCP}.

\begin{table}[ht]
   \centering

   \begin{tabular}{l|c|c|c}
      Method & Residual - $\|A\hat{x}-b\|_2$ & Norm - $\|\hat{x}\|_2$ & Time (s)\\
      \hline
      \texttt{QR}   & $3.0\times10^{-9}$  & $1.3\times10^5$ & $0.04$\\
      \texttt{QRCP} & $2.5\times10^{-13}$ & $5.8\times10^0$ & $0.19$\\
      BLENDENPIK    & $1.4\times10^{-13}$ & $1.4\times10^0$ & $0.16$\\
      \RURVHaar     & $5.8\times10^{-12}$ & $1.5\times10^2$ & $0.52$\\
      \RURVROS      & $1.4\times10^{-12}$ & $4.3\times10^1$ & $0.10$
   \end{tabular}
   \caption{
      Comparison of basic solution residuals for the $1000\times1500$ matrix from Listing \ref{lst:corr_mat} with $p=10$ correlated columns.
      As expected, unpivoted \texttt{QR} has a relatively large residual, while the other methods perform better.
      Note that BLENDENPIK computes the minimum norm solution.
   }
   \label{tab:basic_resids}
\end{table}

\begin{table}[ht]
   \centering
   \begin{tabular}{l|c|c|c}
      Method & Residual - $\|A\hat{x}-b\|_2$ & Norm - $\|\hat{x}\|_2$ & Time (s)\\
      \hline
      \texttt{QR}   & $5.4\times10^{-13}$ & $1.8\times10^1$ & $0.04$\\
      \texttt{QRCP} & $2.5\times10^{-13}$ & $6.1\times10^0$ & $0.20$\\
      BLENDENPIK    & $1.3\times10^{-13}$ & $1.4\times10^0$ & $0.15$\\
      \RURVHaar     & $4.8\times10^{-12}$ & $1.2\times10^2$ & $0.52$\\
      \RURVROS      & $1.3\times10^{-12}$ & $3.9\times10^1$ & $0.10$
   \end{tabular}
   \caption{
      Comparison of basic solution residuals for the $1000\times1500$ matrix from Listing \ref{lst:corr_mat} with $p=0$ correlated columns.
      All methods perform well, and that the two URV-based methods compute mixed basic solutions with norms on the same order as in the previous case with correlated columns.
   }
   \label{tab:basic_resids2}
\end{table}

\section{Experimental Comparison of \RURVHaar~and \RURVROS}
\label{sec:rr}

In this section we experiment with a variety of QR and URV factorizations, some of which are known to be rank-revealing.
In Subsection \ref{ssec:scaling_rr_cond} we experiment with how the rank-revealing conditions \eqref{eq:rr_cond_ratios} and \eqref{eq:rr_cond_strong} scale with increasing $n$.
Our chief interest here is the comparison of \RURVHaar~and \RURVROS.

We can use \RURVROS~to form low-rank approximations by performing the mixing and pre-sort as usual, but only performing $k$ steps of the QR factorization, yielding a rank-$k$ approximation.
The mixing step costs $\mathcal{O}(mn\log n)$ FLOPs as usual, but now the partial QR factorization costs only $\mathcal{O}(mnk)$ FLOPs.
In Subsection \ref{ssec:qlp_approx}, we investigate pairing QR and URV factorizations with Stewart's QLP approximation to the SVD \cite{stewart1999qlp}.
One can use the QLP approximation to obtain an improved rank-$k$ approximation by truncating $L$ factor.

\subsection{Scaling of Rank-Revealing Conditions}
\label{ssec:scaling_rr_cond}

It was shown in \cite{demmel2007fast, ballard2010minimizing} that \RURVHaar~produces a strong rank-revealing factorization with high probability.
\RURVROS~simply replaces Haar mixing with ROS mixing and adds a pre-sort before the unpivoted QR factorization, so we expect \RURVROS~to behave similarly to \RURVHaar.
Specifically, we hope that \RURVROS~obeys the strong rank-revealing conditions \eqref{eq:rr_cond_ratios} and \eqref{eq:rr_cond_strong} in a manner similar to \RURVHaar.

\begin{figure}[ht!]
   \centering
   \includegraphics[width=0.7\textwidth]{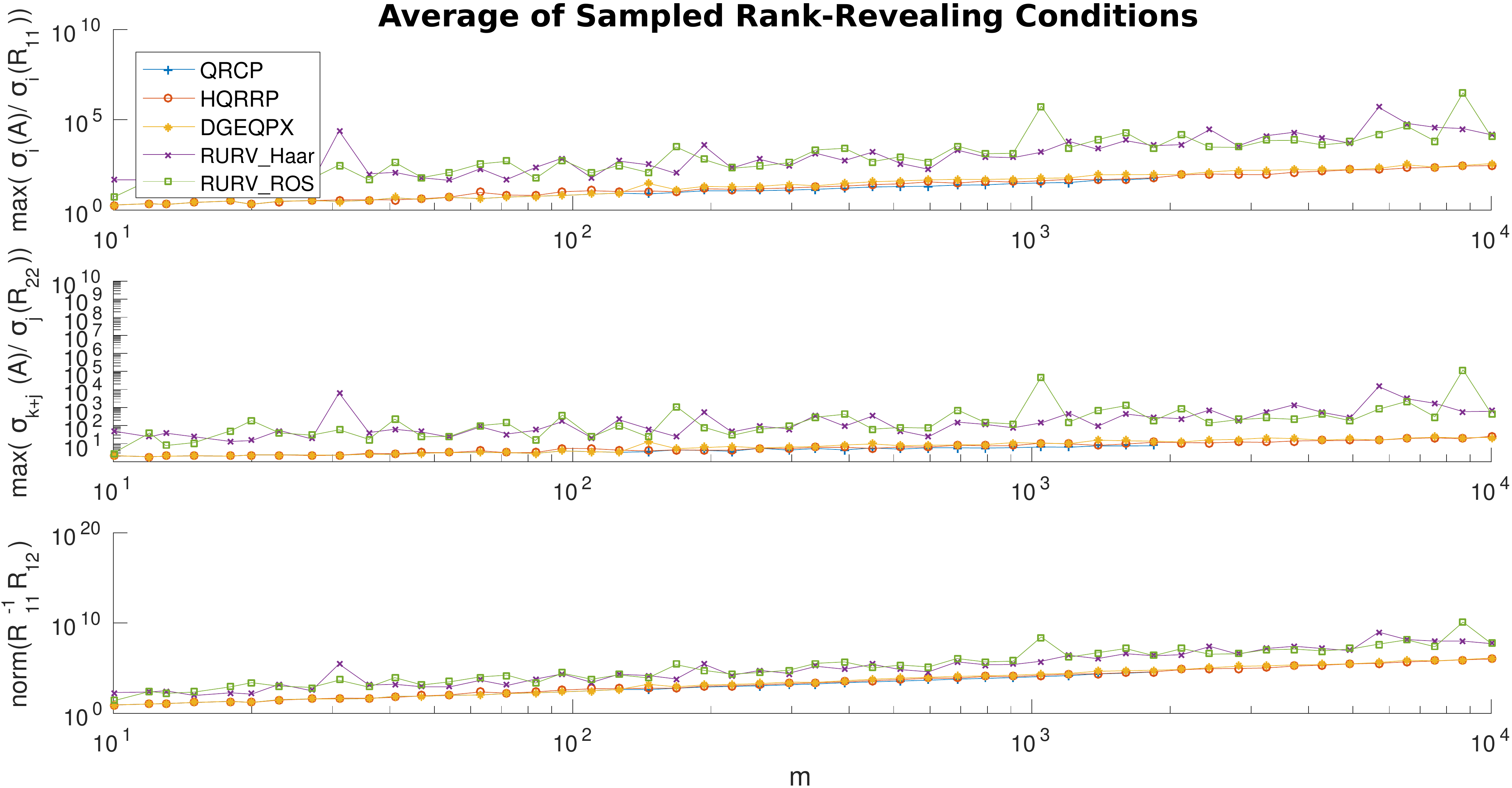}
   \caption{
      Maximum of the sampled values of the rank-revealing conditions from Subsection \ref{ssec:RURVHaar} for five random $m\times m$ matrices of numerical rank $m/2$.
      The three QR factorizations exhibit growth that is clearly bounded by a slowly growing polynomial (linear in a log-log plot).
      \RURVHaar~and \RURVROS~also appear to exhibit bounded growth, with only a few exceptions; recall that \RURVHaar~produces a strong rank-revealing factorization with high probability, not deterministically.
   }
   \label{fig:rr_conditions_scaling_rankk}
\end{figure}

We experimentally test the scaling of the ratios $\sigma_i(A)/\sigma_i(R_{11})$, $\sigma_j(R_{22})/\sigma_{k+j}(A)$, and the norm $\|R_{11}^{-1}R_{12}\|$ to determine if they appear to be bounded above by a slowly-growing polynomial.
In Figure \ref{fig:rr_conditions_scaling_rankk} we take $A$ to be a random $m\times m$ matrix of rank $k\approx m/2$.
The matrix $A$ is formed as $A=U\Sigma V^T$, where $U$ and $V$ are Haar random orthogonal matrices, $\Sigma = \operatorname{diag}(\sigma_1,...,\sigma_m)$, and the $\sigma_i$ decay slowly until $\sigma_{m/2}$, where there is a gap of about $10^{-10}$, after which the $\sigma_i$ decay slowly again.
We sample sizes $m$ from 10 to 1000; for each $m$, we generate five instantiations of the matrix $A$, perform a variety of factorizations for each $A$, and compute the conditions \eqref{eq:rr_cond_ratios} and \eqref{eq:rr_cond_strong} for each factorization.
For plotting, we plot the maximum over the five instantiations of $\max_i\sigma_i(A)/\sigma_i(R_{11})$, $\max_j\sigma_j(R_{22})/\sigma_{k+j}(A)$, and $\|R_{11}^{-1}R_{12}\|$.

We use the highly-accurate LAPACK routine \texttt{dgejsv} to compute singular values of the test matrices (when the exact singular values are unknown) and in the computation of the ratios \eqref{eq:rr_cond_ratios}.
\texttt{dgejsv} implements a preconditioned Jacobi SVD algorithm, which can be more accurate  for small singular values \cite{drmavc2008new, drmavc2008new2}.
Specifically, if $A=DY$ (or $A=YD$), where $D$ is a diagonal weighting matrix and $Y$ is reasonably well-conditioned, \texttt{dgejsv} is guaranteed to achieve high accuracy.
The relative error of the singular values computed with the preconditioned Jacobi method are $\mathcal{O}(\epsilon)\kappa_2(Y)$, whereas the relative errors as computed with a QR-iteration based SVD are $\mathcal{O}(\epsilon)\kappa_2(A)$ \cite{drmavc2008new, demmel2015communication}.
This fact is particularly relevant when we test with the Kahan matrix, which is discussed later in the section.
Even when $A$ is not of the form $A=DY$, $A=YD$, or even $A=D_1YD_2$, it is expected that \texttt{dgejsv} returns singular values at least as accurate as a QR-iteration based SVD.

We test \texttt{QRCP}, \RURVHaar, \RURVROS, \texttt{HQRRP} from \cite{martinsson2015householder}, which uses random projections to select blocks of pivots, and \texttt{DGEQPX} from \cite{bischof1998algorithm}, which is known to be a rank-revealing QR.
Note that \texttt{HQRRP} is intended to cheaply produce a column-pivoted Householder QR; it is not a rank-revealing QR, but it tends to be rank-revealing in practice, like \texttt{QRCP}.

Figure \ref{fig:rr_conditions_scaling_rankk} shows the rank-revealing conditions for $A$ a random $m\times m$ matrix of rank $k\approx m/2$.
The three QR factorizations we test, \texttt{QRCP}, \texttt{HQRRP}, and \texttt{DGEQPX}, perform very well, meaning that the sampled rank-revealing conditions appear to be bounded above by a slowly growing polynomial.
Note that Figure \ref{fig:rr_conditions_scaling_rankk} uses a log-log scale, on which polynomial growth appears linear.
As we expect, \RURVROS~performs about as well as \RURVHaar.
With the exception of a few points, \RURVHaar~and \RURVROS~appear to be bounded above by a slowly growing polynomial, albeit a significantly larger polynomial than for the three QR factorizations.
The exceptions may very well be points where \RURVHaar~or \RURVROS~failed to produce a rank-revealing factorization for at least one of the five sampled $A$ matrices.

\begin{figure}[ht!]
   \centering
   \includegraphics[width=0.7\textwidth]{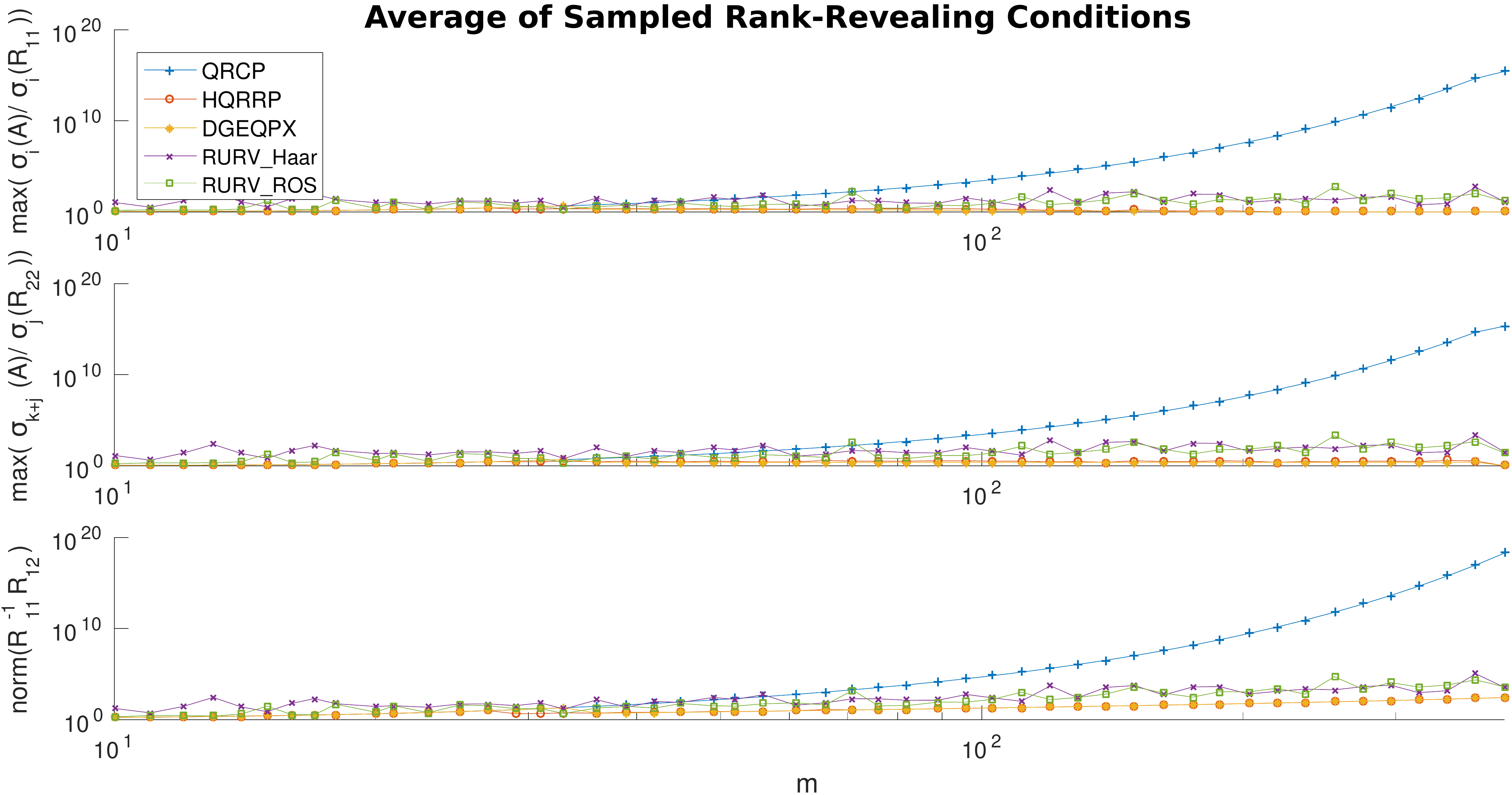}
   \caption{
      Maximum of the sampled values of the rank-revealing conditions from Subsection \ref{ssec:RURVHaar} on the $m\times m$ Kahan matrix.
      As expected, \texttt{QRCP} performs very poorly, with all conditions scaling exponentially.
      Again, we see that \RURVHaar~and \RURVROS~behave similarly.
   }
   \label{fig:rr_conditions_scaling_kahan_rankk}
\end{figure}

Figure \ref{fig:rr_conditions_scaling_kahan_rankk} shows the rank-revealing conditions with $A$ the $m\times m$ Kahan matrix and $k$ chosen to be $m-1$.
The Kahan matrix is a well-known counterexample on which \texttt{QRCP} performs no pivoting in exact arithmetic \cite{drmavc2008failure}.
We use the Kahan matrix (with perturbation) as described in \cite{demmel2015communication}.
The $m\times m$ Kahan matrix is formed as

\begin{equation}
   \label{eq:kahan}
   A = \begin{bmatrix} 1&0&0&\cdots&\cdots&0\\
                       0&s&0&\cdots&\cdots&0\\
                       0&0&s^2&\ddots&\cdots&0\\
                       \vdots&\vdots&\ddots&\ddots&\ddots&\vdots\\
                       \vdots&\vdots&\cdots&\ddots&\ddots&0\\
                       0&0&0&\cdots&0&s^{m-1}
       \end{bmatrix}
       \begin{bmatrix} 1&-c&-c&\cdots&\cdots&-c\\
                       0&1&-c&\cdots&\cdots&-c\\
                       0&0&1&\ddots&\cdots&-c\\
                       \vdots&\vdots&\ddots&\ddots&\ddots&\vdots\\
                       \vdots&\vdots&\cdots&\ddots&\ddots&-c\\
                       0&0&0&\cdots&0&1
       \end{bmatrix},
\end{equation}

\noindent where $s^2+c^2=1$ and $s,c\ge 0$.
When using \texttt{QRCP} to compute the factorization

\[ A\Pi = QR = Q\begin{bmatrix}R_{11}&R_{12}\\&R_{12}\end{bmatrix}, \quad R_{11}\in\reals^{k\times k}, R_{12}\in\reals^{k\times (m-k)}, R_{22}\in\reals^{(m-k)\times(m-k)}, \]

\noindent it is known that $\sigma_k(A)/\sigma_k(R_{11}) \ge \frac{1}{2}c^3(1+c)^{m-4}/s$ for $k=m-1$, and $\sigma_k(R_{11})$ can be much smaller than $\sigma_k(A)$ \cite{gu1996efficient}.
That is, \texttt{QRCP} does not compute a rank-revealing factorization, as the first ratio in \eqref{eq:rr_cond_ratios} grows exponentially for $i=k=m-1$.
To prevent \texttt{QRCP} from pivoting on the Kahan matrix in finite arithmetic, we multiply the $j$th column by $(1-\tau)^{j-1}$, with $1\gg \tau \gg \epsilon$ \cite{drmavc2008failure, demmel2015communication}.
In our tests, we pick $c=0.1$ and $\tau=10^{-7}$.

The most apparent feature of Figure \ref{fig:rr_conditions_scaling_kahan_rankk} is that the rank-revealing conditions for \texttt{QRCP} grow exponentially.
This is a known feature of the Kahan matrix, and shows that \texttt{QRCP} is not \emph{strictly speaking} a rank-revealing QR (in practice, however, it is still used as a rank-revealing factorization).
Moreover, the Kahan matrix is so bad for \texttt{QRCP}, we believe \texttt{dgejsv} cannot accurately compute the singular values in the ratios $\sigma_i(A)/\sigma_i(R_{11})$ and $\sigma_j(R_{22})/\sigma_{k+j}(A)$.
As $m$ grows, the right-hand matrix in \eqref{eq:kahan} becomes increasingly ill-conditioned, and we see the exponential growth in Figure \ref{fig:rr_conditions_scaling_kahan_rankk} stop around $m=10^3$.
In infinite precision arithmetic, the exponential growth should continue, so we stop testing at $m\approx 400$.
As expected, the rank-revealing conditions for \RURVROS~scale in the same manner as \RURVHaar, giving credence to our thought that \RURVROS~is rank-revealing with high probability.

\subsection{Accuracy of R-Values}
Another test we perform involves the accuracy of $|R(i,i)|$ in predicting $\sigma_i(R)$ ($R(i,i)$ is the $i$th diagonal element of the upper-triangular factor from a QR or URV).
Following \cite{demmel2015communication}, we call the $|R(i,i)|$ R-values.
The R-values can be used as a rough estimate of the singular values.
A better approximation is to use Stewart's QLP factorization \cite{stewart1999qlp}, which we discuss in Subsection \ref{ssec:qlp_approx}.
Nevertheless, it is descriptive to investigate the behavior of the R-values.

\begin{figure}[ht!]
   \centering
   \includegraphics[width=0.9\textwidth]{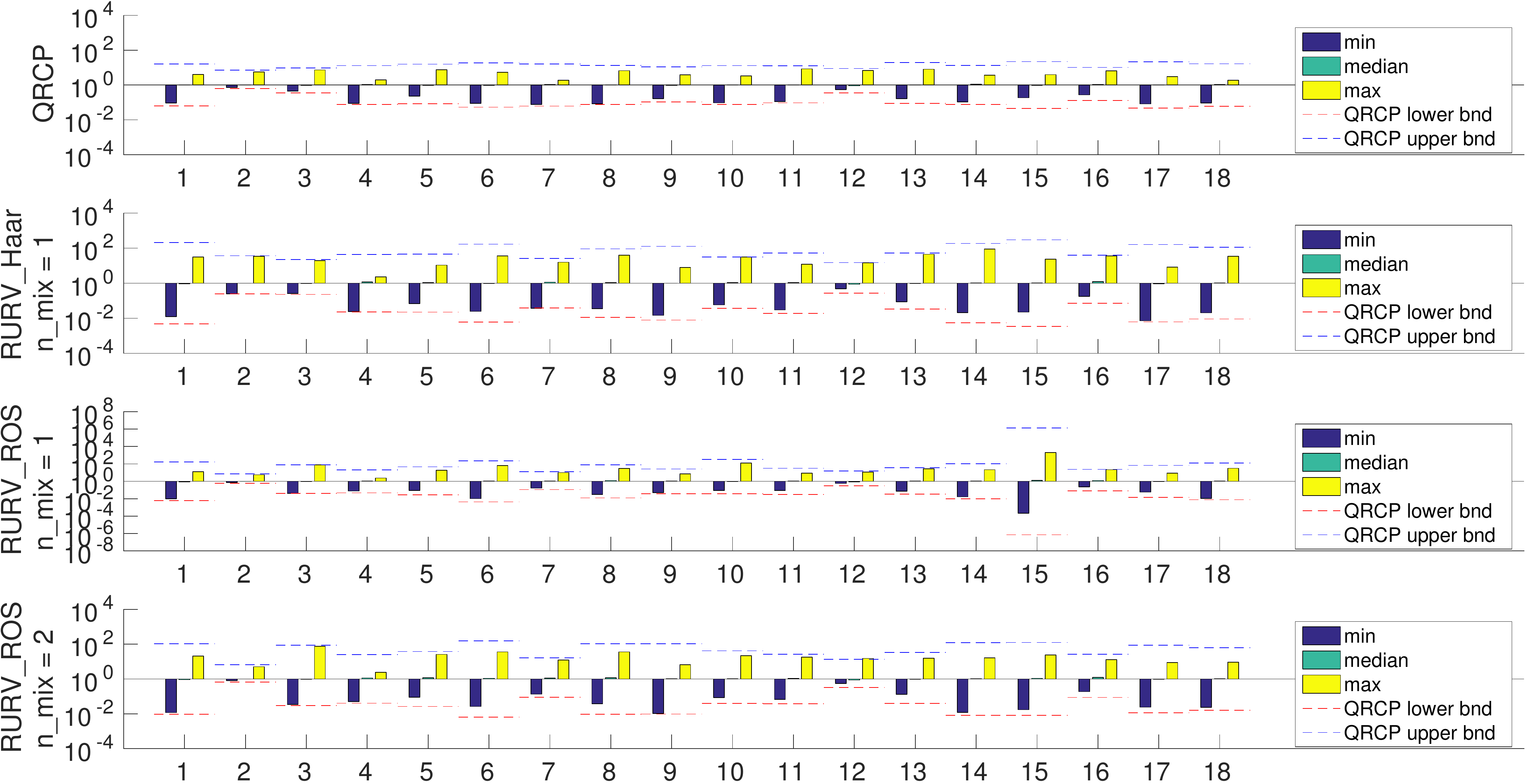}
   \caption{
      Ratios of $|R(i,i)|/\sigma_i(R)$ for the 18 matrices in Table 2 of \cite{demmel2015communication}.
      The abscissa is the index of the matrix in the test set.
      For matrix 15 (SPIKES), \RURVHaar~produces ratios on par with the rest of the test set.
      For \RURVROS, however, using only 1 mixing step produces very bad max/min ratios; using two two mixing steps produces better ratios, but more mixing steps doesn't appear to yield further improvements.
      For the other 17 matrices, \RURVROS~produces ratios comparable with \RURVHaar.
   }
   \label{fig:rvalues_bars}
\end{figure}

We test \texttt{QRCP}, \RURVHaar, and \RURVROS~on the first 18 test matrices from Table 2 of \cite{demmel2015communication} (most matrices are from \cite{hansen2007regularization, gu1996efficient}).
In Figure \ref{fig:rvalues_bars}, we plot the minimum, median, and maximum of the ratios $|R(i,i)|/\sigma_i$ for the 18 test matrices.
For each matrix, we let \texttt{r} and \texttt{s} be the vectors of R-values and singular values, respectively; we plot \texttt{min(r./s)}, \texttt{median(r./s)}, and \texttt{max(r./s)} (using {\sc matlab} syntax).
We see that \texttt{QRCP} produces ratios that are at most just over an order of magnitude away from one.
\RURVHaar~produces slightly worse ratios, which seem to be spread over about two orders of magnitude away from one.
\RURVROS~with one mixing iteration produces ratios comparable to \RURVHaar, with the exception of matrix 15, SPIKES.
For matrix 15, the extreme ratios are significantly larger than on the rest of the test set.
Adding a second mixing iteration brings the ratios back down to a couple orders of magnitude away from one, but does not improve the ratios for the other matrices beyond what is accomplished with a single mixing.  We can also find a bound for the ratios obtained with QR and URV factorizations.

Let $D$ be the diagonal part of $R$ obtained from a QR or URV factorization, and define $Y$ via $R=DY^T$.  This results in the factorization $A\Pi=QDY^T$ for \texttt{QRCP} and $A=UDY^TV$ for \RURVHaar~and \RURVROS.
For \texttt{QRCP}, the diagonal elements of $R$ are non-negative and sorted in decreasing order; this is not guaranteed for \RURVHaar~or \RURVROS.
It follows from the Courant-Fischer minimax theorem \cite{golub1998matrix} that \texttt{QRCP} has the bounds

\begin{equation}
   \label{eq:QRCP_bounds}
   \dfrac{1}{\|Y\|} \le \dfrac{R(i,i)}{\sigma_i} \le \|Y^{-1}\|.
\end{equation}

\noindent For \RURVHaar~and \RURVROS, let $\rho_i$ be the $i$th largest (in absolute value) diagonal element of $R$.
For \RURVHaar~and \RURVROS, we have the bounds

\[ \dfrac{1}{\|Y\|} \le \dfrac{|\rho_i|}{\sigma_i} \le \|Y^{-1}\|. \] 

In addition to the minimum, median, and maximum values of $|R(i,i)|/\sigma_i$ for each matrix, we plot the bounds \eqref{eq:QRCP_bounds} for both \texttt{QRCP} and the two RURV factorizations.
Even though the two RURV factorizations are not guaranteed to be bound by \eqref{eq:QRCP_bounds}, since it is a strong rank-revealing URV, we expect the R-values to somewhat closely approximate the singular values and approximately obey the \texttt{QRCP} bounds.
With the exception of matrix 12, formed as \verb|A=2*rand(n)-1| in {\sc matlab}, we see this behavior in Figure \ref{fig:rvalues_bars}, and we again see \RURVROS~behaving similarly to \RURVHaar.

\subsection{Experiments With the QLP Approximation}
\label{ssec:qlp_approx}
The QLP factorization was introduced by G.W. Stewart as an approximation to the SVD in \cite{stewart1999qlp}.
The idea of the pivoted QLP factorization is to use \texttt{QRCP} to find R-values, and then improve the accuracy (by a surprising amount) by performing another \texttt{QRCP} on $R^T$.
This results in a factorization of the form $A=Q_1\Pi_2LQ_2^T\Pi_1$, where $L$ is lower triangular.
Following \cite{demmel2015communication}, we call the diagonal elements of the $L$ matrix L-values.
In Stewart's original experiments, it was found that L-values approximate the singular values significantly more accurately than the R-values.
Also, the accuracy seemed intimately tied to using \texttt{QRCP} for the first factorization, but that unpivoted \texttt{QR} could be used in the second QR factorization with only cosmetic differences.
It was later shown that the QLP factorization can be interpreted as the first two steps of QR-style SVD algorithm \cite{huckaby2003convergence}.

\begin{figure}[ht!]
   \centering
   \includegraphics[width=0.49\textwidth]{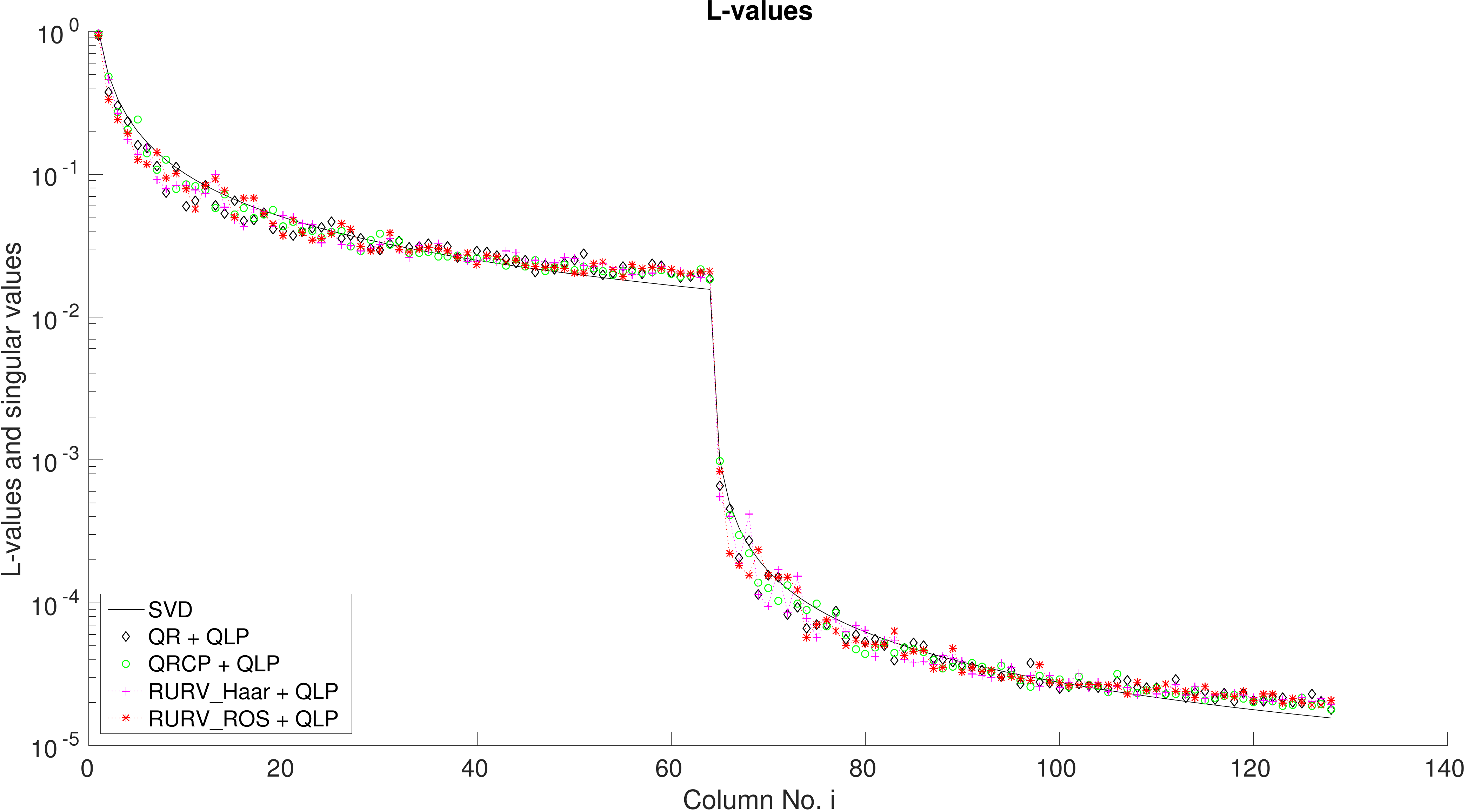}
   \includegraphics[width=0.49\textwidth]{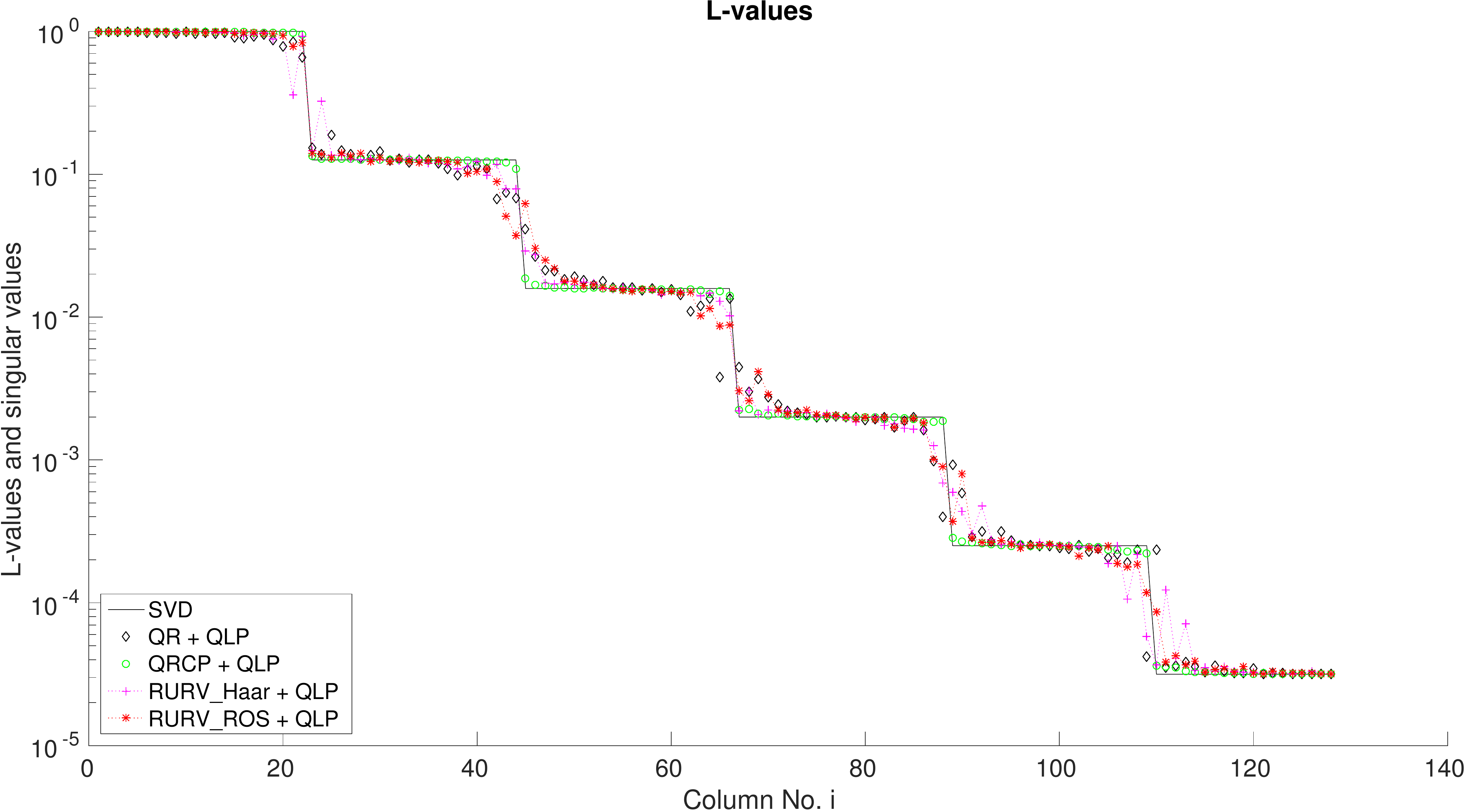}
   \caption{
      L-values from various QLP factorizations on a random $128\times 128$ matrix with slowly decaying singular values and a small gap of approximate size $10^{-1}$ and the Devil's stairs with gaps of approximate size $10^{-1}$.
   In each case, the legend name is of the form \{factorization\}+QLP, where the second factorization is always unpivoted \texttt{QR}.
   }
   \label{fig:lval_rank_devil}
\end{figure}

We experiment with QLP-style factorizations by performing \texttt{QR}, \texttt{QRCP}, \RURVHaar, or \RURVROS, and following up with an unpivoted \texttt{QR} to compute the L-values.
We denote such a factorization as \{factorization\}+QLP (e.g., \texttt{QRCP}+QLP).
For the RURV factorizations, this QLP-style factorization is of the form $A=ULQ^TV$.
Figure \ref{fig:lval_rank_devil} shows the singular values and L-values for a random matrix of the form $A=U\Sigma V^T$, where $U,V$ are Haar random orthogonal, and the singular values are chosen to decay slowly, have a gap of approximate width $10^{-1}$, and decay slowly again.
We see that all QLP-style factorizations, including \texttt{QR}+QLP, identify both the location and magnitude of the gap quite accurately.

Also shown in Figure \ref{fig:lval_rank_devil} are the L-values for the Devil's stairs matrix, which is a particularly difficult example for rank-revealing factorizations.
The Devil's stairs matrix is discussed in \cite{stewart1999qlp, demmel2015communication}, and is formed with $A=U\Sigma V^T$, with $U,V$ Haar random orthogonal and $\Sigma$ controlling the stair-step behavior.
Of all the factorizations, \texttt{QRCP}+QLP performs the best, accurately identifying the location and size of the singular value gaps.
\texttt{QR}+QLP, \RURVHaar+QLP, and \RURVROS+QLP all provide evidence for the existence of singular value gaps, but none is able to identify the precise location and size of the gaps.

\begin{figure}[ht!]
   \centering
   \includegraphics[width=0.48\textwidth]{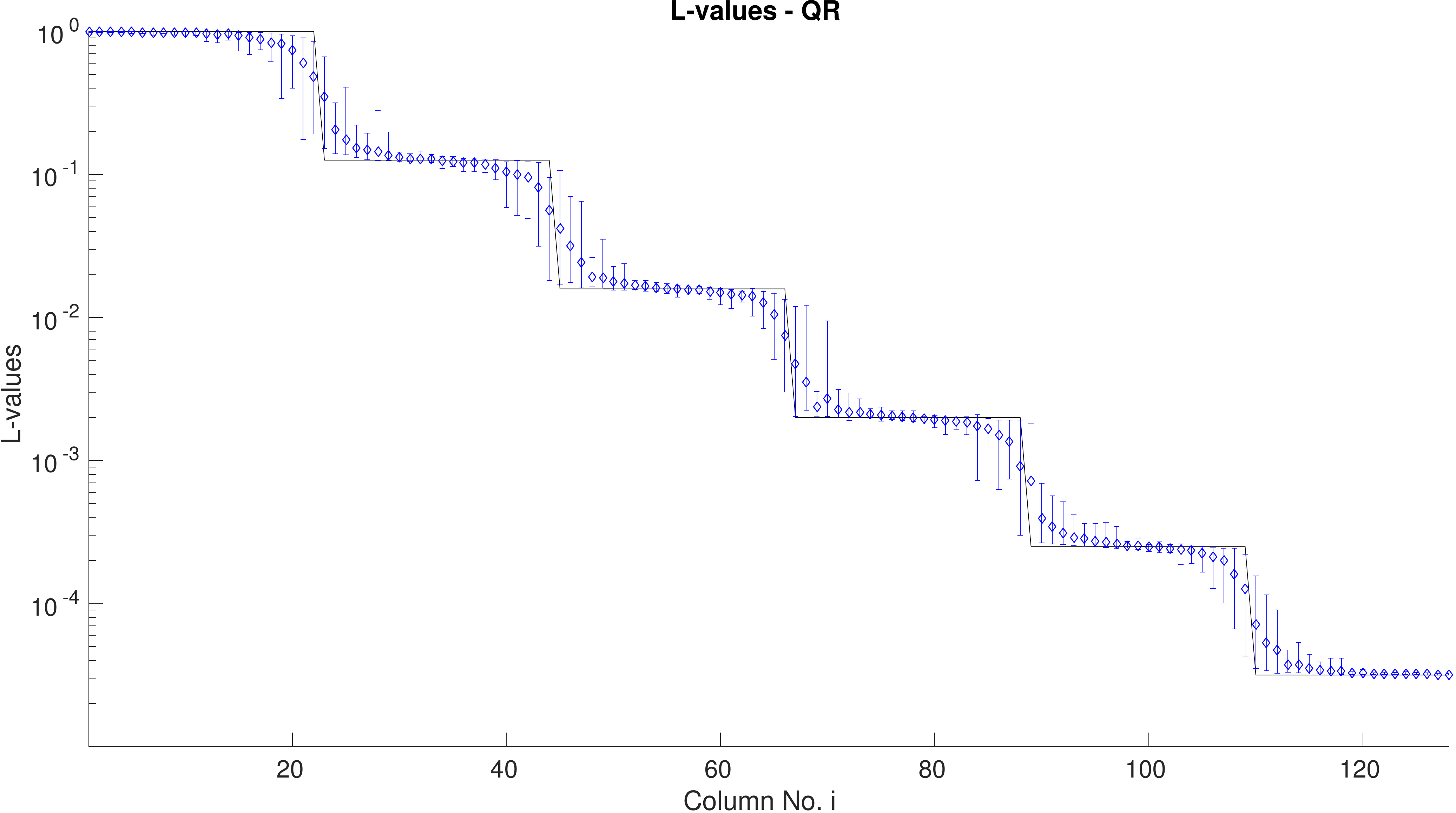}
   \includegraphics[width=0.48\textwidth]{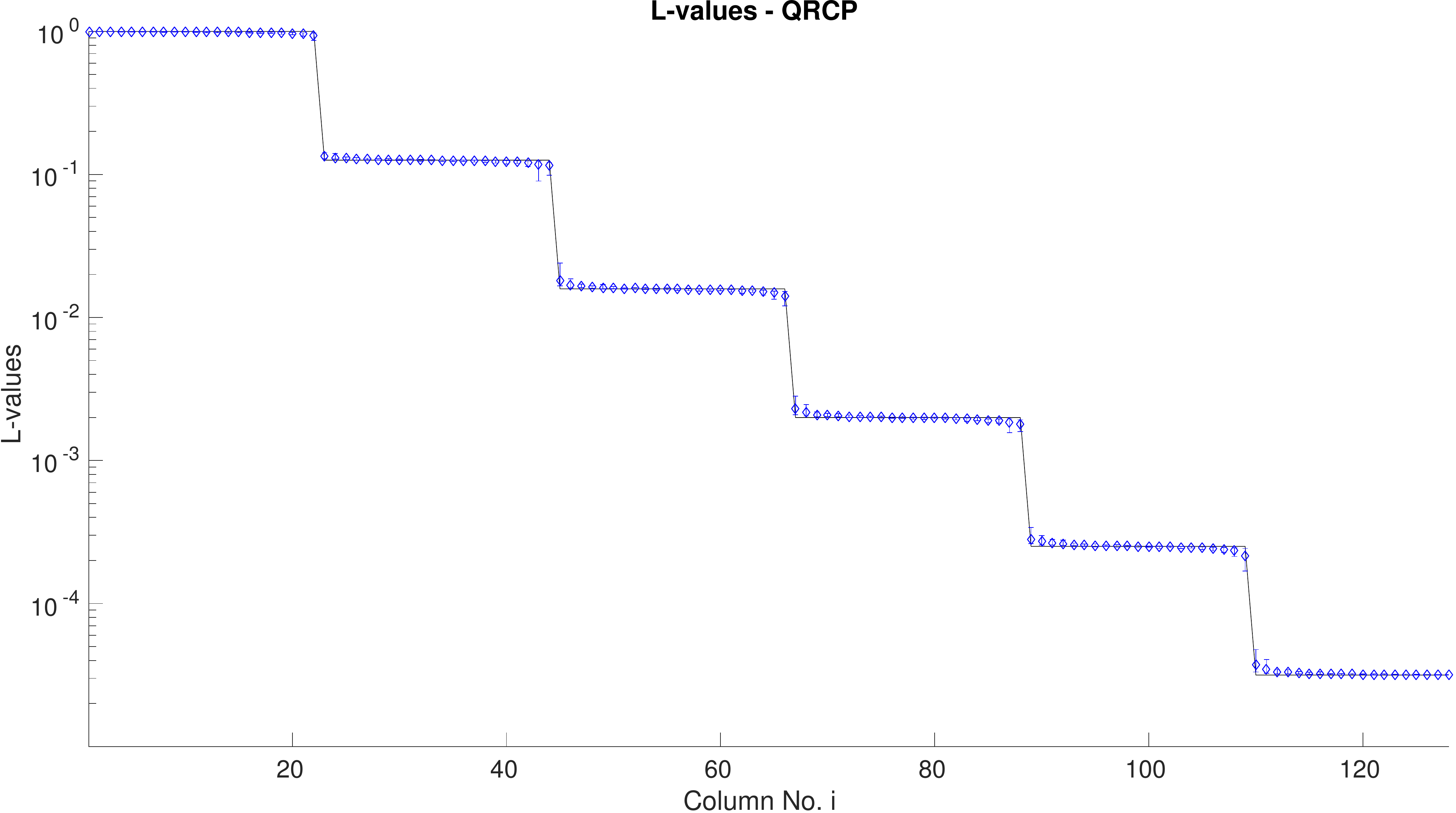}
   \includegraphics[width=0.48\textwidth]{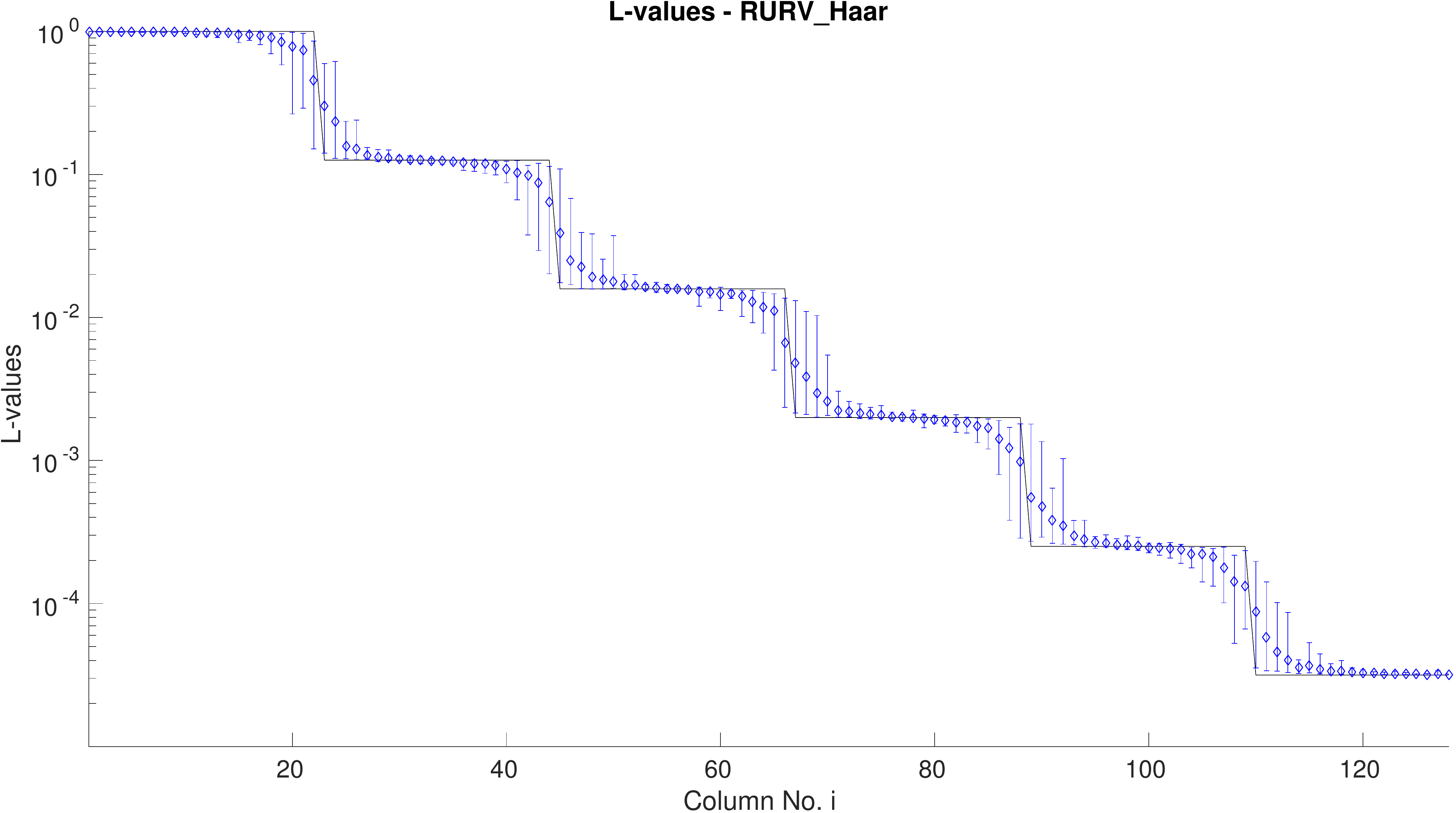}
   \includegraphics[width=0.48\textwidth]{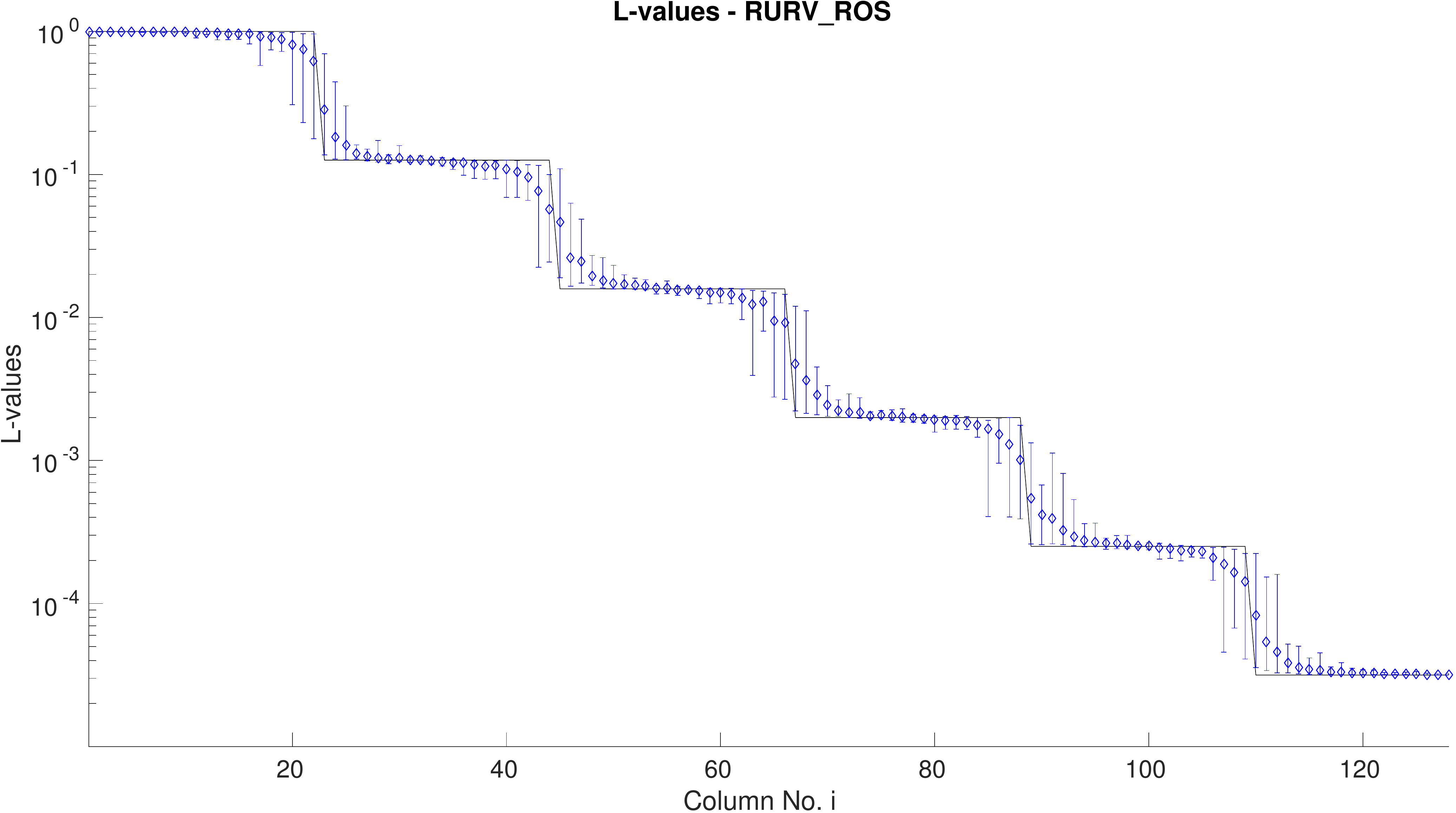}
   \caption{
      Min/Median/Max L-values for 25 runs of a randomly generated $128\times 128$ Devil's stairs matrix with jumps of approximate size $10^{-1}$.
      \texttt{QR}, \RURVHaar, and \RURVROS~appear to predictably show the presence and approximate location of the gaps, but are not accurate enough to estimate the size of the gaps.
      \texttt{QRCP} performs very well, and accurately shows the location and size of the gaps.
   }
   \label{fig:med_lvals}   
\end{figure}

Figure \ref{fig:med_lvals} shows the minimum, median, and maximum L-values for 25 realizations of the Devil's stairs matrix.
We again use \texttt{QR}, \texttt{QRCP}, \RURVHaar, and \RURVROS, followed by \texttt{QR} to form the QLP factorization.
It is clear that \texttt{QRCP}+QLP produces the best L-values;  \RURVHaar+QLP and \RURVROS+QLP generate L-values visually similar to those produced with \texttt{QR}+QLP.
The L-values are smeared around the jumps for \texttt{QR} and the two RURV factorizations, but the L-values have a lower variance around the middle of the flat stairs.
The variance of the L-values around the gaps appears visually similar for \texttt{QR}+QLP, \RURVHaar+QLP, and \RURVROS+QLP.
For \texttt{QR}+QLP, the variance is explained only by the Haar random orthogonal matrices used to construct the Devil's stairs matrix; for \RURVHaar+QLP and \RURVROS+QLP, the variance is a combination of the random Devil's stairs matrix and the random mixing.

\section{Discussion}
\label{sec:discussion}

We have modified \RURVHaar, a strong rank-revealing factorization with high probability, to use random orthogonal mixing (ROS) instead of Haar orthogonal matrix mixing.
The new algorithm, \RURVROS, applies the mixing matrix implicitly and quickly, as opposed to \RURVHaar, where the mixing matrix is generated with an unpivoted \texttt{QR} and applied with dense matrix-matrix multiplication.
With both randomized URV factorizations, one of the principal attractions is the use of cheaper, unpivoted \texttt{QR}, instead of relying on the more expensive \texttt{QRCP}.
The ansatz is that mixing reduces the variance of the column norms, reducing the effect that column pivoting would have, and so we can forgo pivoting and use a cheaper, unpivoted QR.
A URV factorization can be used in many applications that call for a QR, and since the dominant asymptotic cost of \RURVROS~is the same as unpivoted \texttt{QR}, \RURVROS~has the potential to be used as a safer alternative to unpivoted \texttt{QR}.
We have considered only real matrices, but the extension to complex matrices and transforms is natural.

We experiment with using \RURVROS~to solve over- and underdetermined least-squares problems.
Using a URV factorization to solve least-squares is very similar to using a QR factorization.
Our implementation of \RURVROS~even performs comparably to BLENDENPIK, which uses mixing and row sampling to create a preconditioner for LSQR.

When one wants a solution to an underdetermined system, but does not need the minimum norm solution, \RURVROS~can be used to find a basic solution slightly faster than BLENDENPIK, which computes the minimum norm solution.
Additionally, if even a few of the columns of the $A$ matrix are highly correlated, using unpivoted \texttt{QR}, or \texttt{QRCP} on the first $m$ columns, can lead to an inaccurate basic solution; using \RURVROS~computes a mixed basic solution with an accurate residual and for which the norm of the solution is only an order of magnitude larger than the minimum norm solution.

Finally, we experiment with the possible rank-revealing nature of \RURVROS.
We test the scaling of the rank-revealing conditions \eqref{eq:rr_cond_ratios} and \eqref{eq:rr_cond_strong} for \RURVHaar, \RURVROS, and a few other QR factorizations, one of which is rank-revealing.
The prominent feature of the scaling tests is that \RURVROS~behaves very similarly to \RURVHaar, which leads us to suspect that \RURVROS~produces a strong rank-revealing factorization with high probability.
We plan to investigate theoretically the apparent rank-revealing nature of \RURVROS.

\bibliographystyle{alpha}
\newcommand{\etalchar}[1]{$^{#1}$}

\end{document}